\documentclass[10pt]{article}
\usepackage{tikz}
\tikzset{dashdot/.style={dash pattern=on 4pt off 1pt on 1pt off 1pt}}
\usepackage{url}
\usepackage{amsthm,amssymb}
\usepackage{amsfonts}
\usepackage{verbatim}
\usepackage{setspace}
\usepackage{float}
\theoremstyle{plain}
\newtheorem{theorem}{Theorem}
\newtheorem{lemma}[theorem]{Lemma}
\newtheorem{proposition}[theorem]{Proposition}
\newtheorem{corollary}[theorem]{Corollary}
\newtheorem{remark}[theorem]{Remark}
\theoremstyle{definition}
\newtheorem{definition}[theorem]{Definition}

\DeclareSymbolFont{AMSb}{U}{msb}{m}{n}
\DeclareMathSymbol{\N}{\mathbin}{AMSb}{"4E}
\DeclareMathSymbol{\Z}{\mathbin}{AMSb}{"5A}
\DeclareMathSymbol{\R}{\mathbin}{AMSb}{"52}
\DeclareMathSymbol{\Q}{\mathbin}{AMSb}{"51}
\DeclareMathSymbol{\I}{\mathbin}{AMSb}{"49}
\DeclareMathSymbol{\C}{\mathbin}{AMSb}{"43}
\newcommand{\ip}[2]{\ensuremath{\left\langle{#1,#2}\right\rangle}}
 
\title{Cyclotomic Matrices and Graphs over the ring of integers of some imaginary quadratic fields.}
\author{Graeme Taylor}
\date{}
\begin{document}
\maketitle
\begin{abstract} \noindent We determine all Hermitian $\mathcal{O}_{\Q(\sqrt{d})}$-matrices for which every eigenvalue is in the interval $[-2,2]$, for each $d\in \{-2,-7,-11,-15\}$. To do so, we generalise charged signed graphs to $\mathcal{L}$-graphs for appropriate finite sets $\mathcal{L}$, and classify all $\mathcal{L}$-graphs satisfying the same eigenvalue constraints. We find that, as in the integer case, any such matrix / graph is contained in a maximal example with all eigenvalues $\pm 2$.
\end{abstract}

\section{Introduction}
Given a monic polynomial $P(z)=\prod_{i=1}^{d}(z-\alpha_i)\in\Z[z]$, the \emph{Mahler Measure} $M(P)$ is given by \[M(P):=\prod_{i=1}^d \max{(1,|\alpha_i|)}=\prod_{|\alpha_i|>1} |\alpha_i|\] Clearly $M(P)\ge1$; by a result of Kronecker \cite{Kron} $M(P)=1$ if and only if $\pm P$ is the product of a cyclotomic polynomial\footnote{Following Boyd \cite{Boyd}, we will use `cyclotomic' to refer to any polynomial for which all roots are roots of unity, rather than just the irreducible examples $\Phi_n$.} and a power of $z$. For a monic integer polynomial with $M(P)>1$, Lehmer asked (in \cite{Lehmer}) whether $M(P)$ could be arbitrarily close to $1$. This is now known as \emph{Lehmer's Problem}; the negative result - that there is some $\lambda>1$ such that $M(P)>1 \Rightarrow M(P)\ge\lambda$ - is sometimes referred to as \emph{Lehmer's Conjecture}. 

For a monic polynomial $g\in\Z[x]$ of degree $n$, define its \emph{associated reciprocal polynomial} to be $z^ng(z+1/z)$ which is a monic reciprocal polynomial of degree $2n$. For $A$ an $n$-by-$n$ symmetric matrix with entries from $\Z$, denote by $R_A(z)$ the associated reciprocal polynomial of its characteristic polynomial $\chi_A(x)=\mbox{det}(xI-A)$. Further, define $M(A)$, \emph{the Mahler measure of $A$}, to be $M(R_A(z))$. Then $R_A(z)$ has Mahler measure $1$ precisely when $A$ has spectral radius at most 2; we therefore describe such an $A$ as a \emph{cyclotomic} matrix.

McKee and Smyth classified all cyclotomic integer symmetric matrices in \cite{McSm}; in \cite{McSm2} they were then able to prove \begin{equation}\label{lehmerforISM}M(A)\ge\lambda_0=1.17628\ldots\end{equation} for any noncyclotomic integer symmetric matrix $A$. Results of Breusch (\cite{Breusch}) and Smyth (\cite{Smyth2}) prove Lehmer's Conjecture for nonreciprocal monic polynomials with integer coefficients; (\ref{lehmerforISM}) would complete the proof if for every monic reciprocal polynomial $P\in\Z[z]$ there existed an integer symmetric matrix $A$ such that $M(P)=M(A)$. 

Clearly, this would hold if for every $P\in\Z[z]$ there existed an integer symmetric matrix $A$ with $P=R_A(z)$, but counterexamples are easily constructed by identifying polynomials that cannot be the characteristic polynomial of any integer symmetric matrix. In \cite{EsGu} Estes and Guralnick demonstrated that if $f\in\Z[x]$ is a monic, separable, degree $n\le4$ polynomial with all real roots, then $f$ is the minimal polynomial of a $(2n)\times(2n)$ integer symmetric matrix. They thus conjectured that for such $f$ of any degree there is an integer symmetric matrix with $f$ as minimal polynomial. In \cite{Dob} Dobrowolski proves that this is not so, even with the relaxation of the dimension condition: there are infinitely many algebraic integers whose minimal polynomial is not the minimal polynomial of an integer symmetric matrix. 

The results of \cite{McSm2} go further: there it is shown that if an integer symmetric matrix $A$ is noncyclotomic with $M(A)<1.3$, then $M(A)$ is one of sixteen given values. By comparison with the tables of small Salem numbers (\cite{Boyd}, \cite{Moss}), noncyclotomic counterexamples to the existence of an $A$ satisfying $M(A)=M(P)$ for any given $P$ are found: the polynomial $z^{14}-z^{12}+z^7-z^2+1$ has $M(P)=1.20261\ldots$, but this is not one of the possible $M(A)<1.3$ if $A$ is an integer symmetric matrix. 

Lehmer's problem therefore remains open for reciprocal polynomials due to these `missing' Mahler measures. An obvious approach is to extend the study of integer symmetric matrices to broader classes of combinatorial objects that still yield integer polynomials. In this paper we take the first step in extending to Hermitian matrices with entries from the rings of integers of various imaginary quadratic fields, by classifying all cyclotomic examples over these rings also. 

As shall be seen, it suffices to classify all maximal connected cyclotomic $\mathcal{L}$-graphs for appropriate finite sets $\mathcal{L}$. Theorems \ref{cyc_class_neg2}, \ref{cyc_class_neg7}, \ref{cyc_class_neg11} and \ref{cyc_class_neg15} present such a classification for $\mathcal{O}_{\Q(\sqrt{-2})}$, $\mathcal{O}_{\Q(\sqrt{-7})}$, $\mathcal{O}_{\Q(\sqrt{-11})}$ and $\mathcal{O}_{\Q(\sqrt{-15})}$ respectively.

\section{Cyclotomic Integer Symmetric Matrices}

If $A$ is a block diagonal matrix, then its list of eigenvalues is the union of the lists of the eigenvalues of the blocks. If there is a reordering of the rows (and columns) of $A$ such that it has block diagonal form with more than one block, then $A$ will be called \emph{decomposable}; if there is no such reordering, $A$ is called \emph{indecomposable}. Clearly any decomposable cyclotomic matrix decomposes into cyclotomic blocks, so to classify all cyclotomic matrices it is sufficent to identify the indecomposable ones. 

The following result is of central importance to this effort:

\begin{theorem}[Cauchy Interlacing Theorem\footnote{This is Th\'eor\`eme I of Cauchy's curiously titled paper \cite{Cauchy} from 1829. For a modern reference in English see Theorem 4.3.8 of \cite{Horn}, which provides a proof by the Courant-Fischer min-max theorem (\emph{Id.} Theorem 4.2.11); a very short proof by reduction to interlacing of polynomials is given in Fisk \cite{Fisk}.} ]\label{cauchyinterlacetheorem}
Let $A$ be a Hermitian $n\times n$ matrix with eigenvalues $\lambda_1\le\lambda_2\le\cdots\le\lambda_n$. \\Let $B$ be obtained from $A$ by deleting row $i$ and column $i$ from $A$. \\Then the eigenvalues $\mu_1\le\cdots\le\mu_{n-1}$ of $B$ interlace with those of $A$: that is, \[ \lambda_1\le\mu_1\le\lambda_2\le\mu_2\le\cdots\le\lambda_{n-1}\le\mu_{n-1}\le\lambda_n .\]
\end{theorem}

Thus if $A$ is cyclotomic, so is any $B$ obtained by successively deleting a series of rows and corresponding columns from $A$. We describe such a $B$ as being \emph{contained} in $A$. If an indecomposable cyclotomic matrix $A$ is not contained in a strictly larger indecomposable cyclotomic matrix, then we call $A$ \emph{maximal}.

Additionally, an equivalence relation on cyclotomic matrices can be defined as follows. Let $O_n(\Z)$ denote the orthogonal group of $n\times n$ signed permutation matrices. Conjugation of a cyclotomic matrix by a matrix from this group gives another matrix with the same eigenvalues, which is thus also cyclotomic. Cyclotomic matrices $A,A'$ related in this way are described as \emph{strongly equivalent}; indecomposable cyclotomic matrices $A$ and $A'$ are then considered equivalent if $A'$ is strongly equivalent to either $A$ or $-A$.

The following is an easy consequence of Theorem \ref{cauchyinterlacetheorem}: 

\begin{lemma}[\cite{McSm}, Lemma 6]\label{ISMisCSG} Apart from matrices equivalent to either $(2)$ or $\left(\begin{array}{cc} 0 & 2 \\ 2 & 0\end{array}\right)$, any indecomposable cyclotomic matrix has all entries from the set $\{0,1,-1\}$.
\end{lemma}

This motivates the following generalisations of the adjacency matrix of a graph. If $A$ is an $n\times n$ matrix with diagonal entries all zero and off-diagonal elements from $\{0,1,-1\}$ then $A$ describes an $n$-vertex \emph{signed} graph (as in \cite{CST}, \cite{Zas}), whereby a non-zero $(i,j)$th entry indicates an edge between vertices $i$ and $j$ with a `sign' of $-1$ or $1$.  For a general $\{0,1,-1\}$ matrix we extend this to \emph{charged signed graphs}, interpreting a non-zero diagonal entry as a `charge' on the corresponding vertex. 

A charged signed graph $G$ is therefore described as cyclotomic if its adjacency matrix $A$ is cyclotomic; the Mahler measure of $G$ is that of $A$ (i.e., of $R_A(z)$), and graphs $G,G'$ are (strongly) equivalent if and only if their adjacency matrices $A,A'$ are. A charged signed graph $G$ is connected if and only if its adjacency matrix is indecomposable. If a cyclotomic matrix $A'$ is contained in $A$ then its corresponding charged signed graph $G'$ is an induced subgraph of $G$ corresponding to $A$; thus a maximal cyclotomic charged signed graph is not an induced subgraph of any strictly larger connected cyclotomic charged signed graph.

The equivalence relation on matrices has the following interpretation for graphs. $O_n(\Z)$ is generated by matrices of the form $\mbox{diag}(1,1,\ldots,1,-1,1,\ldots,1)$ and permutation matrices. Conjugation by the former has the effect of negating the signs of all edges incident at some vertex $v$; following \cite{CST} this is described as \emph{switching at $v$}. Conjugation by a permutation matrix merely permutes vertex labels and so up to equivalence we may ignore vertex labellings: strong equivalence classes are therefore determined only by switching operations on unlabelled graphs. Equivalence of charged signed graphs is then generated by switching and the operation of negating all edge signs and vertex charges of a connected component.

For conciseness, we indicate edge signs visually, with a sign of $1$ given by an unbroken line $\begin{tikzpicture}
\node (a) at (0,0) {};
\node (b) at (1,0) {};
\draw (a) -- (b);
\end{tikzpicture}$ and a sign of $-1$ given by a dotted line $\begin{tikzpicture}
\node (a) at (0,0) {};
\node (b) at (1,0) {};
\draw [dotted,thick] (a) -- (b);
\end{tikzpicture}$. Vertices with charge $0$ (neutral), $1$ (positive) and $-1$ (negative) will be drawn as $\begin{tikzpicture}
\node (a) at (0,0) [fill=black,draw,shape=circle] {};
\end{tikzpicture}$,
$\begin{tikzpicture}
\node (a) at (0,0) [fill=white,inner sep=0pt,draw,shape=circle] {$+$};
\end{tikzpicture}$ and
$\begin{tikzpicture}
\node (a) at (0,0) [fill=white,inner sep=0pt,draw,shape=circle] {$-$};
\end{tikzpicture}$ respectively.

By Lemma \ref{ISMisCSG} we thus have that (with the exception of the given matrices) any maximal indecomposable cyclotomic integer symmetric matrix is the adjacency matrix of a maximal connected cyclotomic charged signed graph.

\section{Maximal Connected Cyclotomic Charged Signed Graphs}
A complete classification of cyclotomic matrices over $\Z$ is therefore given via the main results of \cite{McSm}:

\begin{theorem}[\cite{McSm} Theorem 1]\label{McSmThm1} Every maximal connected cyclotomic signed graph is equivalent to one of the following:
\begin{itemize}
\item[(i)] The 14-vertex signed graph $S_{14}$ shown in \cite{McSm} Fig. 3;
\item[(ii)] The 16-vertex signed graph $S_{16}$ shown in \cite{McSm} Fig. 4;
\item[(iii)] For some $k=3,4,\ldots$, the $2k$-vertex toral tessellation $T_{2k}$ shown in \cite{McSm} Fig. 1.
\end{itemize}
Further, every connected cyclotomic signed graph is contained in a maximal one.
\end{theorem}

\begin{theorem}[\cite{McSm} Theorem 2]\label{McSmThm2} Every maximal connected cyclotomic charged signed graph not included in Theorem \ref{McSmThm1} is equivalent to one of the following:
\begin{itemize}
\item[(i)] One of the three sporadic charged signed graphs $S_7,S_8,S_8'$ shown in \cite{McSm} Fig. 7;
\item[(ii)] For some $k=2,3,4,\ldots$, one of the two $2k$-vertex cylindrical tessellations $C_{2k}^{++},C_{2k}^{+-}$ shown in \cite{McSm} Fig. 6.
\end{itemize}
Further, every connected cyclotomic charged signed graph is contained in a maximal one.
\end{theorem}

\section{Cyclotomic $\mathcal{L}$-graphs}

If we now let $A$ be a Hermitian matrix with all entries from $R=\mathcal{O}_{\Q(\sqrt{d})}$ for $d<0$, then $\chi_A(x)\in\Z[x]$ and so $R_A(z)\in\Z[z]$. Further, Theorem \ref{cauchyinterlacetheorem} still applies, with the following corollary:

\begin{lemma}\label{inL} Let $A$ be an $n\times n$ cyclotomic Hermitian matrix. Then \[|A_{i,j}\overline{A_{i,j}}|\le4\] for all $1\le i,j \le n$.\end{lemma}
\begin{proof} By interlacing, if $(A^2)_{i,i}>4$ for any $i$ then $A^2$ has an eigenvalue $\lambda$ such that $|\lambda|>4$ and thus $A$ has an eigenvalue $\lambda'$ such that $|\lambda'|>2$. Therefore for $A$ to be cyclotomic we require $(A^2)_{i,i}\le4$, which implies \[ |A_{i,j}\overline{A_{i,j}}| \le \sum_{k=1}^{n} A_{i,k}\overline{A_{i,k}} = \sum_{k=1}^{n} A_{i,k}A_{k,i} = (A^2)_{i,i} \le 4\] \end{proof}

For $R=\mathcal{O}_{\Q(\sqrt{d})}$ and $n\ge1$, define $\mathcal{L}_n=\{x\in R\,|\,x\overline{x}=n\}$. If $x=a+b\sqrt{d} \in R$ then $x\overline{x}=a^2-db^2=Norm(x)\in \Z$, so $x=0$ or $x\in\mathcal{L}_n$ for some $n$. Thus if $A$ is a cyclotomic Hermitian matrix with all entries from $R$, then by Lemma \ref{inL} $A$ is an $\mathcal{L}$-matrix for \[ \mathcal{L}:= \{0\} \cup \mathcal{L}_1 \cup \mathcal{L}_2 \cup \mathcal{L}_3 \cup \mathcal{L}_4.\]

\begin{corollary} If $d$ is squarefree and $d\not\in\{-1,-2,-3,-7,-11,-15\}$ then $\mathcal{L}=\{0,\pm1,\pm2\}$ and thus any cyclotomic Hermitian $\mathcal{L}$-matrix is an integer symmetric matrix. \end{corollary}

We restrict our attention to $d$ satisfying $\mathcal{L}$ finite, $\mathcal{L}\neq\{0,\pm1,\pm2\}$ and $\mathcal{L}_1=\{\pm1\}$: that is, $d\in\{-2,-7,-11,-15\}$. The remaining cases ($d=-1,-3$) will be presented in future work. 

As in the $\Z$-matrix case, for $n>1$ we have that an indecomposable cyclotomic $\mathcal{L}$-matrix has diagonal entries from $\{0,1,-1\}$. We may therefore generalise the study of charged signed graphs to charged $\mathcal{L}$-graphs by identifying diagonal entries with charges in the usual way, whilst for $i<j$ a non-zero $(i,j)$th entry $x\in\mathcal{L}$ corresponds to an edge with label $x$ between vertices $i$ and $j$. We inherit the notions of indecomposability and maximality; strong equivalence holds as before, although we also consider all of $A,-A,\overline{A},-\overline{A}$ to be equivalent. 

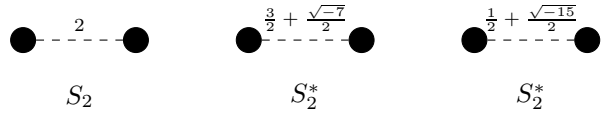
\begin{figure}[H]
\[\begin{tikzpicture}[scale=1.5]
\node (a) at (0,0) [fill=black,draw,shape=circle] {};
\node (b) at (1,0) [fill=black,draw,shape=circle] {};
\draw [dashed] (a) -- node[above] {\scriptsize $2$} (b);
\node (a) at (2,0) [fill=black,draw,shape=circle] {};
\node (b) at (3,0) [fill=black,draw,shape=circle] {};
\draw [dashed] (a) -- node[above] {\scriptsize $\frac{3}{2}+\frac{\sqrt{-7}}{2}$} (b);
\node (a) at (4,0) [fill=black,draw,shape=circle] {};
\node (b) at (5,0) [fill=black,draw,shape=circle] {};
\draw [dashed] (a) -- node[above] {\scriptsize $\frac{1}{2}+\frac{\sqrt{-15}}{2}$} (b);
\node (S2label) at (0.5,-0.5) {$S_2$};
\node (S2label1) at (2.5,-0.5) {$S_2^*$};
\node (S2*label2) at (4.5,-0.5) {$S_2^*$};
\end{tikzpicture}\]
\caption{The 2-vertex sporadic maximal connected cyclotomic $\mathcal{L}$-graphs $S_2$, $S_2^*$.}
\label{weight4figure}
\end{figure}

\begin{figure}[H]
\[\begin{tikzpicture}[scale=1.5]
\node (a) at (0,0) [fill=white,inner sep=0pt,draw,circle,label=left:{\scriptsize$1$}] {$+$};
\node (b) at (1,0) [fill=white,inner sep=0pt,draw,circle,label=right:{\scriptsize$2$}] {$-$};
\draw (a.0) -- (b.180) (a.340) -- (b.200) (a.20) -- node[above] {\scriptsize$1+\sqrt{-2}$} (b.160);
\node (S2label) at (0.5,-0.5) {$S_2'$};
\node (a) at (2,0) [fill=white,inner sep=0pt,draw,circle,label=left:{\scriptsize$1$}] {$+$};
\node (b) at (3,0) [fill=white,inner sep=0pt,draw,circle,label=right:{\scriptsize$2$}] {$-$};
\draw (a.0) -- (b.180) (a.340) -- (b.200) (a.20) -- node[above] {\scriptsize$\frac{1}{2}+\frac{\sqrt{-11}}{2}$} (b.160);
\node (S2label2) at (2.5,-0.5) {$S_2'$};
\end{tikzpicture}
\]
\[\begin{tikzpicture}[scale=1.5]
\node (a) at (0,0) [fill=black,draw,circle,label=above left:{\scriptsize$1$}] {};
\node (b) at (1,0) [fill=black,draw,circle,label=above right:{\scriptsize$2$}] {};
\draw (a.0) -- (b.180) (a.340) -- (b.200) (a.20) -- node[above] {\scriptsize$1+\sqrt{-2}$} (b.160);
\node (c) at (0,-1) [fill=black,draw,circle,label=below left:{\scriptsize$3$}] {};
\node (d) at (1,-1) [fill=black,draw,circle,label=below right:{\scriptsize$4$}] {};
\draw (c.0) -- (d.180) (c.340) -- node[below] {\scriptsize$-1-\sqrt{-2}$} (d.200) (c.20) -- (d.160);
\draw (a) -- (c) (b) -- (d);
\node (S4label) at (0.5,-0.5) {$S_4'$};
\node (a) at (2,0) [fill=black,draw,circle,label=above left:{\scriptsize$1$}] {};
\node (b) at (3,0) [fill=black,draw,circle,label=above right:{\scriptsize$2$}] {};
\draw (a.0) -- (b.180) (a.340) -- (b.200) (a.20) -- node[above] {\scriptsize$\frac{1}{2}+\frac{\sqrt{-11}}{2}$} (b.160);
\node (c) at (2,-1) [fill=black,draw,circle,label=below left:{\scriptsize$3$}] {};
\node (d) at (3,-1) [fill=black,draw,circle,label=below right:{\scriptsize$4$}] {};
\draw (c.0) -- (d.180) (c.340) -- node[below] {\scriptsize$-\frac{1}{2}-\frac{\sqrt{-11}}{2}$} (d.200) (c.20) -- (d.160);
\draw (a) -- (c) (b) -- (d);
\node (S4label2) at (2.5,-0.5) {$S_4'$};
\end{tikzpicture}
\]
\caption{The 2-vertex and 4-vertex sporadic maximal connected cyclotomic charged $\mathcal{L}$-graphs $S_2'$ and $S_4'$.}
\label{weight3figure}
\end{figure}
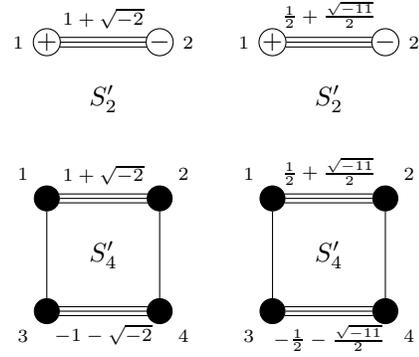

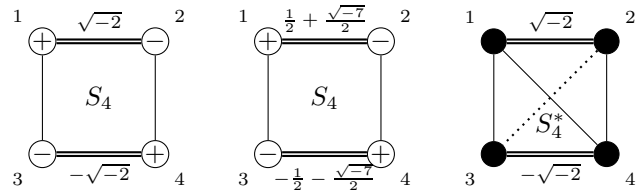
\begin{figure}[H]
\[\begin{tikzpicture}[scale=1.5]
\node (a) at (0,0) [fill=white,inner sep=0pt,draw,circle,label=above left:{\scriptsize$1$}] {$+$};
\node (b) at (1,0) [fill=white,inner sep=0pt,draw,circle,label=above right:{\scriptsize$2$}] {$-$};
\node (c) at (0,-1) [fill=white,inner sep=0pt,draw,circle,label=below left:{\scriptsize$3$}] {$-$};
\node (d) at (1,-1) [fill=white,inner sep=0pt,draw,circle,label=below right:{\scriptsize$4$}] {$+$};
\draw [double,thick] (a) -- node[above] {\scriptsize$\sqrt{-2}$} (b) (c) -- node[below] {\scriptsize$-\sqrt{-2}$} (d);
\draw (a) -- (c) (b) -- (d);
\node (S4label) at (0.5,-0.5) {$S_4$};
\node (a) at (2,0) [fill=white,inner sep=0pt,draw,circle,label=above left:{\scriptsize$1$}] {$+$};
\node (b) at (3,0) [fill=white,inner sep=0pt,draw,circle,label=above right:{\scriptsize$2$}] {$-$};
\node (c) at (2,-1) [fill=white,inner sep=0pt,draw,circle,label=below left:{\scriptsize$3$}] {$-$};
\node (d) at (3,-1) [fill=white,inner sep=0pt,draw,circle,label=below right:{\scriptsize$4$}] {$+$};
\draw [double,thick] (a) -- node[above] {\scriptsize$\frac{1}{2}+\frac{\sqrt{-7}}{2}$} (b) (c) -- node[below] {\scriptsize$-\frac{1}{2}-\frac{\sqrt{-7}}{2}$} (d);
\draw (a) -- (c) (b) -- (d);
\node (S4label) at (2.5,-0.5) {$S_4$};
\node (a) at (4,0) [fill=black,draw,circle,label=above left:{\scriptsize$1$}] {};
\node (b) at (5,0) [fill=black,draw,circle,label=above right:{\scriptsize$2$}] {};
\node (c) at (4,-1) [fill=black,draw,circle,label=below left:{\scriptsize$3$}] {};
\node (d) at (5,-1) [fill=black,draw,circle,label=below right:{\scriptsize$4$}] {};
\draw [double,thick] (a) -- node[above] {\scriptsize$\sqrt{-2}$} (b);
\draw [double,thick] (c) -- node[below] {\scriptsize$-\sqrt{-2}$} (d);
\draw (a) -- (c) (b) -- (d) (a) -- (d);
\draw [dotted,thick] (b) -- (c);
\node (S4label) at (4.5,-0.75) {$S_4^*$};
\end{tikzpicture}
\]
\caption{The 4-vertex sporadic maximal connected cyclotomic charged $\mathcal{L}$-graphs $S_4$ and $S_4^*$.}
\label{S4w2figure}
\end{figure}

We will extend the results of \cite{McSm} to the following:

\begin{theorem}{$(d=-2)$}\label{cyc_class_neg2} Every maximal connected cyclotomic $\mathcal{L}$-graph for $R=\mathcal{O}_{\Q(\sqrt{-2})}$ not included in Theorems \ref{McSmThm1}, \ref{McSmThm2} is equivalent to one of the following:
\begin{itemize}
\item[(i)] The 2-vertex $\mathcal{L}$-graph $S_2$ shown in Fig. \ref{weight4figure} or $S_2'$ shown in Fig. \ref{weight3figure} ;
\item[(ii)] One of the 4-vertex $\mathcal{L}$-graphs $S_4'$, $S_4$ or $S_4^*$ shown in Figs. \ref{weight3figure} and \ref{S4w2figure};
\item[(iii)] The 8-vertex $\mathcal{L}$-graph $S_8^*$ shown in Fig. \ref{S8*figure};
\item[(iv)] For some $k=2,3,4,\ldots$, the $2k$-vertex $\mathcal{L}$-graph $T_{2k}^4$ shown in Fig. \ref{T2k4figure};
\item[(v)] For some $k=1,2,3,\ldots$, the $2k+1$-vertex $\mathcal{L}$-graph $C_{2k}^{2+}$ shown in Fig. \ref{C2k+figure}.
\end{itemize}
\end{theorem}

\begin{theorem}{$(d=-7)$}\label{cyc_class_neg7} Every maximal connected cyclotomic $\mathcal{L}$-graph for $R=\mathcal{O}_{\Q(\sqrt{-7})}$ not included in Theorems \ref{McSmThm1}, \ref{McSmThm2} is equivalent to one of the following:
\begin{itemize}
\item[(i)] The 2-vertex $\mathcal{L}$-graph $S_2$ or $S_2^*$ shown in Fig. \ref{weight4figure};
\item[(iii)] The 4-vertex $\mathcal{L}$-graph $S_4$ shown in Fig. \ref{S4w2figure};
\item[(iv)] The 6-vertex $\mathcal{L}$-graph $S_6^\dag$ shown in Fig. \ref{S6dagfigure};
\item[(v)] The 8-vertex $\mathcal{L}$-graph $S_8^*$ shown in Fig. \ref{S8*figure};
\item[(vi)] For some $k=2,3,4,\ldots$, the $2k$-vertex $\mathcal{L}$-graph $T_{2k}^4$ shown in Fig. \ref{T2k4figure};
\item[(vii)] For some $k=2,3,4,\ldots$, the $2k$-vertex $\mathcal{L}$-graph $T_{2k}^{4'}$ shown in Fig. \ref{T2k4'figure};
\item[(viii)] For some $k=1,2,3,\ldots$, the $2k+1$-vertex $\mathcal{L}$-graph $C_{2k}^{2+}$ shown in Fig. \ref{C2k+figure}.
\end{itemize}
\end{theorem}

\begin{theorem}{$(d=-11)$}\label{cyc_class_neg11} Every maximal connected cyclotomic $\mathcal{L}$-graph for $R=\mathcal{O}_{\Q(\sqrt{-11})}$ not included in Theorems \ref{McSmThm1}, \ref{McSmThm2} is equivalent to one of the following:
\begin{itemize}
\item[(i)] The 2-vertex $\mathcal{L}$-graph $S_2$ shown in Fig. \ref{weight4figure} or $S_2'$ shown in Fig. \ref{weight3figure} ;
\item[(ii)] The 4-vertex $\mathcal{L}$-graph $S_4'$ shown in Fig. \ref{weight3figure}.
\end{itemize}
\end{theorem}

\begin{theorem}{$(d=-15)$}\label{cyc_class_neg15} Every maximal connected cyclotomic $\mathcal{L}$-graph for $R=\mathcal{O}_{\Q(\sqrt{-15})}$ not included in Theorems \ref{McSmThm1}, \ref{McSmThm2} is equivalent to either the 2-vertex $\mathcal{L}$-graph $S_2$ or the 2-vertex $\mathcal{L}$-graph $S_2^*$ as shown in Fig. \ref{weight4figure}.
\end{theorem}

\begin{theorem}\label{containedinmax} Every connected cyclotomic $\mathcal{L}$-graph for $R=\mathcal{O}_{\Q(\sqrt{d})}$, $d\in\{-2,-7,-11,-15\}$ is contained in a maximal one.\end{theorem}

\begin{figure}[H]
\[
\begin{tikzpicture}[scale=1.5]
\node (a) at (0:1) [fill=black,draw,circle,label=right:{\scriptsize$4$}] {};
\node (b) at (60:1) [fill=black,draw,circle,label=above right:{\scriptsize$3$}] {};
\node (c) at (120:1) [fill=black,draw,circle,label=above left:{\scriptsize$2$}] {};
\node (d) at (180:1) [fill=black,draw,circle,label=left:{\scriptsize$1$}] {};
\node (e) at (240:1) [fill=black,draw,circle,label=below left:{\scriptsize$6$}] {};
\node (f) at (300:1) [fill=black,draw,circle,label=below right:{\scriptsize$5$}] {};
\draw (d) -- (c) (a) -- (b) (e) -- (f) (b) -- (e) ;
\draw [dotted,thick] (f) -- (c) (d) -- (a);
\draw [double,thick] (d) -- node[below left] {\scriptsize$\overline{\omega}$} (e) (f) -- node[below right] {\scriptsize$-\omega$} (a);
\draw [double,thick] (b) -- node[above] {\scriptsize$\omega$} (c);
\end{tikzpicture}
\]
\caption{The 6-vertex sporadic maximal connected cyclotomic $\mathcal{L}$-graph $S_6^\dag$. ($\omega=\frac{1}{2}+\frac{\sqrt{-7}}{2}$)}
\label{S6dagfigure}
\end{figure}
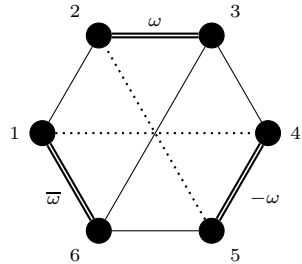

\begin{figure}[H]
\[
\begin{tikzpicture}[scale=2]
\node (1) at (0,0) [fill=black,draw,shape=circle,label=below left:{\scriptsize$8$}] {};
\node (2) at (1,0) [fill=black,draw,shape=circle,label=below right:{\scriptsize$7$}] {};
\node (3) at (0,1) [fill=black,draw,shape=circle,label=below left:{\scriptsize$4$}] {};
\node (4) at (1,1) [fill=black,draw,shape=circle,label=below right:{\scriptsize$3$}] {};
\node (1b) at (0.3,0.3) [fill=black,draw,shape=circle,label=above right:{\scriptsize$5$}] {};
\node (2b) at (1.3,0.3) [fill=black,draw,shape=circle,label=above right:{\scriptsize$6$}] {};
\node (3b) at (0.3,1.3) [fill=black,draw,shape=circle,label=above left:{\scriptsize$1$}] {};
\node (4b) at (1.3,1.3) [fill=black,draw,shape=circle,label=above right:{\scriptsize$2$}] {};
\draw [dotted,thick] (1) -- (2) (1b) -- (2b);
\draw (1) -- (3) (2) -- (4) (1b) -- (3b) (2b) -- (4b) (3) -- (4) (3b) -- (4b);
\draw [double,thick] (1) -- node[right] {\scriptsize$-\sqrt{-2}$} (1b) (2) -- node[right,yshift=-0.5em] {\scriptsize$\sqrt{-2}$} (2b) (4) -- node[left] {\scriptsize$-\sqrt{-2}$} (4b) (3) -- node[left,yshift=0.5em] {\scriptsize$\sqrt{-2}$} (3b);

\node (1) at (2,0) [fill=black,draw,shape=circle,label=below left:{\scriptsize$8$}] {};
\node (2) at (3,0) [fill=black,draw,shape=circle,label=below right:{\scriptsize$7$}] {};
\node (3) at (2,1) [fill=black,draw,shape=circle,label=below left:{\scriptsize$4$}] {};
\node (4) at (3,1) [fill=black,draw,shape=circle,label=below right:{\scriptsize$3$}] {};
\node (1b) at (2.3,0.3) [fill=black,draw,shape=circle,label=above right:{\scriptsize$5$}] {};
\node (2b) at (3.3,0.3) [fill=black,draw,shape=circle,label=above right:{\scriptsize$6$}] {};
\node (3b) at (2.3,1.3) [fill=black,draw,shape=circle,label=above left:{\scriptsize$1$}] {};
\node (4b) at (3.3,1.3) [fill=black,draw,shape=circle,label=above right:{\scriptsize$2$}] {};
\draw [dotted,thick] (1) -- (2) (1b) -- (2b);
\draw (1) -- (3) (2) -- (4) (1b) -- (3b) (2b) -- (4b) (3) -- (4) (3b) -- (4b);
\draw [double,thick] (1) -- node[right] {\scriptsize$-\omega$} (1b) (2) -- node[right,yshift=-0.5em] {\scriptsize$\omega$} (2b) (4) -- node[left] {\scriptsize$-\omega$} (4b) (3) -- node[left,yshift=0.5em] {\scriptsize$\omega$} (3b);
\node (S8label) at (2.5,-0.3) {\scriptsize($\omega=\frac{1}{2}+\frac{\sqrt{-7}}{2}$)};
\end{tikzpicture}
\]
\caption{The 8-vertex sporadic maximal connected cyclotomic $\mathcal{L}$-graphs $S_8^*$.}
\label{S8*figure}
\end{figure}
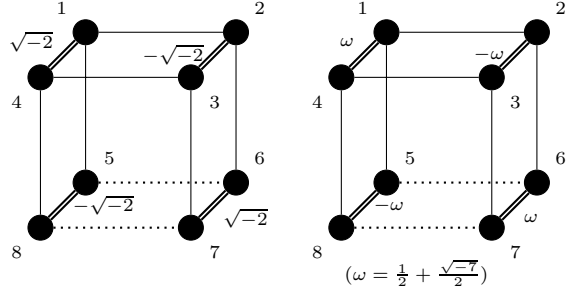

\begin{figure}[H]
\[
\begin{tikzpicture}[scale=1.5]
\node (1) at (0,0) [fill=black,draw,circle,label=above:{\scriptsize$1$}] {};
\node (2) at (1,0) [fill=black,draw,circle,label=above:{\scriptsize$2$}] {};
\node (3) at (2,0) [fill=black,draw,circle,label=above:{\scriptsize$3$}] {};
\node (3a) at (2.4,0) {};
\node (3b) at (2.4,-0.4) {};
\draw (1) -- (2) -- (3) -- (3a); 
\draw (3) -- (3b);
\node (L+1) at (0,-1) [fill=black,draw,circle,label=below:{\scriptsize$k$}] {};
\node (L+2) at (1,-1) [fill=black,draw,circle,label=below:{\scriptsize$k+1$}] {};
\node (L+3) at (2,-1) [fill=black,draw,circle,label=below:{\scriptsize$k+2$}] {};
\node (k+3a) at (2.4,-1) {};
\node (k+3b) at (2.4,-0.6) {};
\draw [dotted] (L+1) -- (L+2) -- (L+3) -- (k+3a);
\draw [dotted] (L+3) -- (k+3b);
\draw (1) -- (L+2) (2) -- (L+3);
\draw [dotted] (L+1) -- (2) (L+2) -- (3);
\node (spacer) at (2.5,-0.5) {$\cdots$};
\node (k-1b) at (2.6,-0.4) {};
\node (2k-1b) at (2.6,-0.6) {};
\node (k-1a) at (2.6,0) {};
\node (L-1) at (3,0) [fill=black,draw,circle,label=above:{\scriptsize$k-2$}] {};
\node (L) at (4,0) [fill=black,draw,circle,label=above:{\scriptsize$k-1$}] {};
\draw (k-1a) -- (L-1) -- (L);
\node (2k-1a) at (2.6,-1) {};
\node (2L-1) at (3,-1) [fill=black,draw,circle,label=below:{\scriptsize$2k-2$}] {};
\node (2L) at (4,-1) [fill=black,draw,circle,label=below:{\scriptsize$2k-1$}] {};
\draw [dotted] (2k-1a) -- (2L-1) -- (2L);
\draw (2k-1b) -- (2L-1) (L-1) -- (2L);
\draw [dotted] (k-1b) -- (L-1) (2L-1) -- (L); 
\node (2L+1) at (-0.5,-0.5) [fill=black,draw,circle,label=left:{\scriptsize$2k-1$}] {};
\node (2L+2) at (4.5,-0.5) [fill=black,draw,circle,label=right:{\scriptsize$2k$}] {};
\draw [double,thick] (L) -- node[black,above right] {\scriptsize $\omega$} (2L+2) -- node[black,below right] {\scriptsize $-\omega$} (2L);
\draw [double,thick] (1) -- node[black,above left] {\scriptsize $\omega$} (2L+1) -- node[black,below left] {\scriptsize $\omega$} (L+1);
\end{tikzpicture}
\]
\caption{The family $T_{2k}^4$ of $2k$-vertex maximal connected cyclotomic $\mathcal{L}$-graphs. \emph{($k\ge2$; $\omega=\sqrt{-2},\frac{1}{2}+\frac{\sqrt{-7}}{2}$ for $d=-2,-7$ respectively.)}}
\label{T2k4figure}
\end{figure}

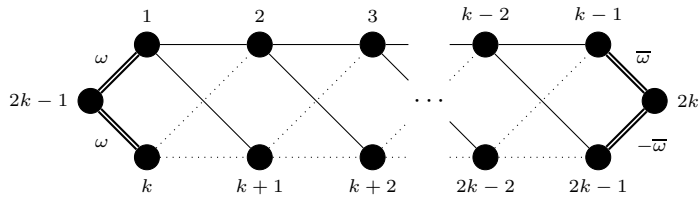
\begin{figure}[H]
\[
\begin{tikzpicture}[scale=1.5]
\node (1) at (0,0) [fill=black,draw,circle,label=above:{\scriptsize$1$}] {};
\node (2) at (1,0) [fill=black,draw,circle,label=above:{\scriptsize$2$}] {};
\node (3) at (2,0) [fill=black,draw,circle,label=above:{\scriptsize$3$}] {};
\node (3a) at (2.4,0) {};
\node (3b) at (2.4,-0.4) {};
\draw (1) -- (2) -- (3) -- (3a); 
\draw (3) -- (3b);
\node (L+1) at (0,-1) [fill=black,draw,circle,label=below:{\scriptsize$k$}] {};
\node (L+2) at (1,-1) [fill=black,draw,circle,label=below:{\scriptsize$k+1$}] {};
\node (L+3) at (2,-1) [fill=black,draw,circle,label=below:{\scriptsize$k+2$}] {};
\node (k+3a) at (2.4,-1) {};
\node (k+3b) at (2.4,-0.6) {};
\draw [dotted] (L+1) -- (L+2) -- (L+3) -- (k+3a);
\draw [dotted] (L+3) -- (k+3b);
\draw (1) -- (L+2) (2) -- (L+3);
\draw [dotted] (L+1) -- (2) (L+2) -- (3);
\node (spacer) at (2.5,-0.5) {$\cdots$};
\node (k-1b) at (2.6,-0.4) {};
\node (2k-1b) at (2.6,-0.6) {};
\node (k-1a) at (2.6,0) {};
\node (L-1) at (3,0) [fill=black,draw,circle,label=above:{\scriptsize$k-2$}] {};
\node (L) at (4,0) [fill=black,draw,circle,label=above:{\scriptsize$k-1$}] {};
\draw (k-1a) -- (L-1) -- (L);
\node (2k-1a) at (2.6,-1) {};
\node (2L-1) at (3,-1) [fill=black,draw,circle,label=below:{\scriptsize$2k-2$}] {};
\node (2L) at (4,-1) [fill=black,draw,circle,label=below:{\scriptsize$2k-1$}] {};
\draw [dotted] (2k-1a) -- (2L-1) -- (2L);
\draw (2k-1b) -- (2L-1) (L-1) -- (2L);
\draw [dotted] (k-1b) -- (L-1) (2L-1) -- (L); 
\node (2L+1) at (-0.5,-0.5) [fill=black,draw,circle,label=left:{\scriptsize$2k-1$}] {};
\node (2L+2) at (4.5,-0.5) [fill=black,draw,circle,label=right:{\scriptsize$2k$}] {};
\draw [double,thick] (L) -- node[black,above right] {\scriptsize $\overline{\omega}$} (2L+2) -- node[black,below right] {\scriptsize $-\overline{\omega}$} (2L);
\draw [double,thick] (1) -- node[black,above left] {\scriptsize $\omega$} (2L+1) -- node[black,below left] {\scriptsize $\omega$} (L+1);
\end{tikzpicture}
\]
\caption{The family $T_{2k}^{4'}$ of $2k$-vertex maximal connected cyclotomic $\mathcal{L}$-graphs. \emph{($k\ge2$, $\omega=\frac{1}{2}+\frac{\sqrt{-7}}{2}$.)}}
\label{T2k4'figure}
\end{figure}

\begin{figure}[H]
\[
\begin{tikzpicture}[scale=1.5]
\node (1) at (0,0) [fill=white,inner sep=0pt,draw,circle,label=above:{\scriptsize$1$}] {$+$};
\node (2) at (1,0) [fill=black,draw,circle,label=above:{\scriptsize$2$}] {};
\node (3) at (2,0) [fill=black,draw,circle,label=above:{\scriptsize$3$}] {};
\node (3a) at (2.4,0) {};
\node (3b) at (2.4,-0.4) {};
\draw (1) -- (2) -- (3) -- (3a); 
\draw (3) -- (3b);
\node (L+1) at (0,-1) [fill=white,inner sep=0pt,draw,circle,label=below:{\scriptsize$k+1$}] {$+$};
\node (L+2) at (1,-1) [fill=black,draw,circle,label=below:{\scriptsize$k+2$}] {};
\node (L+3) at (2,-1) [fill=black,draw,circle,label=below:{\scriptsize$k+3$}] {};
\draw (1) -- (L+1);
\node (k+3a) at (2.4,-1) {};
\node (k+3b) at (2.4,-0.6) {};
\draw [dotted] (L+1) -- (L+2) -- (L+3) -- (k+3a);
\draw [dotted] (L+3) -- (k+3b);
\draw (1) -- (L+2) (2) -- (L+3);
\draw [dotted] (L+1) -- (2) (L+2) -- (3);
\node (spacer) at (2.5,-0.5) {$\cdots$};
\node (k-1b) at (2.6,-0.4) {};
\node (2k-1b) at (2.6,-0.6) {};
\node (k-1a) at (2.6,0) {};
\node (L-1) at (3,0) [fill=black,draw,circle,label=above:{\scriptsize$k-1$}] {};
\node (L) at (4,0) [fill=black,draw,circle,label=above:{\scriptsize$k$}] {};
\draw (k-1a) -- (L-1) -- (L);
\node (2k-1a) at (2.6,-1) {};
\node (2L-1) at (3,-1) [fill=black,draw,circle,label=below:{\scriptsize$2k-1$}] {};
\node (2L) at (4,-1) [fill=black,draw,circle,label=below:{\scriptsize$2k$}] {};
\draw [dotted] (2k-1a) -- (2L-1) -- (2L);
\draw (2k-1b) -- (2L-1) (L-1) -- (2L);
\draw [dotted] (k-1b) -- (L-1) (2L-1) -- (L); 
\node (2L+1) at (4.5,-0.5) [fill=black,draw,circle,label=right:{\scriptsize$2k+1$}] {};
\draw [double,thick] (L) -- node[black,above right] {\scriptsize $\omega$} (2L+1) -- node[black,below right] {\scriptsize $-\omega$} (2L);
\end{tikzpicture}
\]
\caption{The family $C_{2k}^{2+}$ of $2k+1$-vertex maximal connected cyclotomic $\mathcal{L}$-graphs. \emph{($k\ge1$; $\omega=\sqrt{-2},\frac{1}{2}+\frac{\sqrt{-7}}{2}$ for $d=-2,-7$ respectively.)}}
\label{C2k+figure}
\end{figure}

\section{Sporadic $\mathcal{L}$-graphs}

\subsection{Growing cyclotomic $\mathcal{L}$-graphs}
\begin{definition}\label{edgeweight} For an edge with label $x$ we define its \emph{weight} to be the norm of $x$ (so a \emph{weight n edge} is one with a label from $\mathcal{L}_n$). For a vertex $v$, we define its \emph{weighted degree} as the sum of the weights of the edges incident at $v$, plus 1 if $v$ has a charge of $\pm1$. \end{definition}

\begin{proposition}\label{maxdeg4} If $v$ is a vertex in a cyclotomic $\mathcal{L}$-graph, then $v$ has weighted degree at most 4.\end{proposition}

We will often specify the edges of an $\mathcal{L}$-graph only up to their weight; we describe such a representation as the \emph{form} of the graph. Edges without an explicit label will be indicated by dashes ($\begin{tikzpicture}
\node (a) at (0,0) [fill=black,draw,shape=circle] {};
\node (b) at (0.9,0) [fill=black,draw,shape=circle] {};
\draw [dashed] (a) -- (b);
\end{tikzpicture}$, $\begin{tikzpicture}
\node (a) at (0,0) [fill=black,draw,shape=circle] {};
\node (b) at (0.9,0) [fill=black,draw,shape=circle] {};
\draw [dashed,double,thick] (a) -- (b);
\end{tikzpicture}$, $\begin{tikzpicture}
\node (a) at (0,0) [fill=black,draw,shape=circle] {};
\node (b) at (0.9,0) [fill=black,draw,shape=circle] {};
\draw [dashed] (a.0) -- (b.180) (a.15) -- (b.165) (a.345) -- (b.195);
\end{tikzpicture}$ for edges from $\mathcal{L}_1, \mathcal{L}_2, \mathcal{L}_3$ respectively) whilst an unspecified - possibly absent - edge will be shown as $\begin{tikzpicture}
\node (a) at (0,0) [fill=black,draw,shape=circle] {};
\node (b) at (1,0) [fill=black,draw,shape=circle] {};
\draw [dashdot] (a) -- (b); \end{tikzpicture}$. If a vertex is of unknown charge $c\in\{0,1,-1\}$ then we denote it by $\circledast$; a vertex known to be charged but of unknown polarity is denoted {\scriptsize\textcircled{$\pm$}}.

Given an induced subgraph $H$ of a cyclotomic $\mathcal{L}$-graph $G$, we can recover $G$ by reintroducing each missing vertex. By interlacing, each graph in this sequence is itself cyclotomic. Theoretically, any cyclotomic $\mathcal{L}$-graph can therefore be grown from the seed set of 2-vertex $\mathcal{L}$-graphs. The combinatorial explosion in possible vertex additions renders this infeasible as a fully general approach. But we are able to first eliminate higher weight edges from consideration, then with refinement identify induced subgraphs that yield only finitely many maximal cyclotomic $\mathcal{L}$-graphs. Such refinements include reducing modulo equivalence after each round (whilst feasible); ignoring additions that would necessarily yield noncyclotomic examples by Proposition \ref{maxdeg4}; and reducing the search space by fixing edges via switching both in $H$ and the added vertices, which for any $G$ inducing a subgraph of form $H$ will ensure we recover some $G'$ equivalent to $G$.  

By the choice of $d$, if $G$ is an $\mathcal{L}$-graph with all edge labels from $\mathcal{L}_1$ then it is a charged signed graph as classified in \cite{McSm}. Thus we may assume that $G$ has at least one edge label from $\mathcal{L}_2\cup\mathcal{L}_3\cup\mathcal{L}_4$. 

\subsection{$\mathcal{L}$-Graphs with edge labels from $\mathcal{L}_3\cup\mathcal{L}_4$}\label{w34edgeclassification}

\subsubsection{Edge labels from $\mathcal{L}_4$}\label{L4section}
By Proposition \ref{maxdeg4}, if vertices $u,v$ are joined by an edge of weight 4, then they can have no other neighbours. Thus a maximal connected $\mathcal{L}$-graph with a weight 4 edge is necessarily of the form \[\begin{tikzpicture}
\node (a) at (0,0) [fill=black,draw,shape=circle] {};
\node (b) at (1,0) [fill=black,draw,shape=circle] {};
\draw [dashdot] (a) -- node[above] {\scriptsize $t$} (b);
 \end{tikzpicture}\]
where $t\in\mathcal{L}_4$. For $d=-2,-11$ $\mathcal{L}_4=\{\pm2\}$, so such a graph is equivalent to $S_2$ as given in Fig. \ref{weight4figure}. For $d=-7,-15$ we have $\mathcal{L}_4=\{\pm2,\pm 3/2\pm\sqrt{-7}/2\}$ , $\{\pm2,\pm 1/2\pm\sqrt{-15}/2\}$ respectively; up to equivalence if $t\neq\pm2$ then we may assume it is as given for the graphs $S_2^*$ in Fig. \ref{weight4figure}. 

We may therefore restrict our attention to $\mathcal{L}=\mathcal{L}_3\cup\mathcal{L}_2\cup\mathcal{L}_1\cup\{0\}$. Moreover, this completes the classification for $d=-15$, where $\mathcal{L}_2=\mathcal{L}_3=\emptyset$, so Theorem \ref{cyc_class_neg15} holds.

\subsubsection{Edge labels from $\mathcal{L}_3$}\label{L3section}

Let $G$ be a maximal connected cyclotomic $\mathcal{L}$-graph with a weight 3 edge label. For $d=-2$ or $-11$, we have (by negating and/or conjugating if necessary) that $G$ is equivalent to such a graph with an edge label of $\alpha=1+\sqrt{-2}$ or $\alpha=1/2+\sqrt{-11}/2$ respectively. We may thus take as seed set representatives of the cyclotomic graphs of the form $\begin{tikzpicture}
\node (b) at (1,0) [fill=white,inner sep=0pt,draw,shape=circle] {$*$};
\node (a) at (0,0) [fill=white,inner sep=0pt,draw,shape=circle] {$*$};
\draw (a.20) -- node[above] {\scriptsize $\alpha$}  (b.160);
\draw (a.0) -- (b.180);
\draw (a.340) -- (b.200);
\end{tikzpicture}$. The growing algorithm terminates after three rounds, indicating that there are only finitely many maximal cyclotomic $\mathcal{L}$-graphs with a weight 3 edge label. Up to form, they are either $
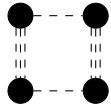
\begin{tikzpicture}
\node (a) at (0,0) [fill=white,inner sep=0pt,draw,shape=circle] {$\pm$};
\node (b) at (1,0) [fill=white,inner sep=0pt,draw,shape=circle] {$\pm$};
\draw [dashed] (a.20) -- (b.160);
\draw [dashed] (a.0) -- (b.180);
\draw [dashed] (a.340) -- (b.200);
\end{tikzpicture}$
or 
\[\begin{tikzpicture}
\node (a) at (2,0) [fill=black,draw,shape=circle] {};
\node (b) at (3,0) [fill=black,draw,shape=circle] {};
\node (c) at (2,-1) [fill=black,draw,shape=circle] {};
\node (d) at (3,-1) [fill=black,draw,shape=circle] {};
\draw [dashed] (a) -- (b) (c) -- (d);
\draw [dashed] (a.250) -- (c.110);
\draw [dashed] (a.270) -- (c.90);
\draw [dashed] (a.290) -- (c.70);
\draw [dashed] (b.250) -- (d.110);
\draw [dashed] (b.270) -- (d.90);
\draw [dashed] (b.290) -- (d.70);
\end{tikzpicture}
\]

It is then straightforward to determine equivalence class representatives; any cyclotomic $\mathcal{L}$-graph of one of the above forms is equivalent to either $S_2'$ or $S_4'$ as given in Fig. \ref{weight3figure}. 

We may therefore restrict our attention to $\mathcal{L}=\mathcal{L}_2\cup\mathcal{L}_1\cup\{0\}$. Moreover, this completes the classification for $d=-11$, where $\mathcal{L}_2=\emptyset$, so Theorem \ref{cyc_class_neg11} holds.

\subsection{Sporadic $\mathcal{L}$-Graphs with edge labels from $\mathcal{L}_2$}

\subsubsection{Charge Isolation}

\begin{lemma}\label{maxis4cyc_A} If $G$ is a maximal cyclotomic $\mathcal{L}$-graph inducing a subgraph of the form $\begin{tikzpicture}
\node (a) at (0,0) [fill=white,inner sep=0pt,draw,shape=circle] {$\pm$};
\node (b) at (1,0) [fill=white,inner sep=1pt,draw,shape=circle] {$*$};
\draw [dashed,double,thick] (a) -- (b);
\end{tikzpicture}$ then $G$ is equivalent to either $C_2^{2+}$ (the 3-vertex case of $C_{2k}^{2+}$ given in Fig. \ref{C2k+figure}) or $S_4$ as given in Fig. \ref{S4w2figure}. \end{lemma}
\begin{proof} Growing from representatives of the seed set of cyclotomic $\mathcal{L}$-graphs of form $\begin{tikzpicture}
\node (a) at (0,0) [fill=white,inner sep=0pt,draw,shape=circle] {$\pm$};
\node (b) at (1,0) [fill=white,inner sep=1pt,draw,shape=circle] {$*$};
\draw [dashed,double,thick] (a) -- (b);
\end{tikzpicture}$ terminates after two rounds, with all maximal examples being of claimed form. Testing then confirms that in each case all cyclotomic examples are equivalent to the given representative.\end{proof}

\subsubsection{Non-cyclotomic structures}

\begin{lemma}\label{maxis4cyc_D} There are no cyclotomic $\mathcal{L}$-graphs of the form $\begin{tikzpicture}
\node (a) at (0,0) [fill=black,draw,circle] {};
\node (b) at (1,0) [fill=black,draw,circle] {};
\node (c) at (2,0) [fill=black,draw,circle] {};
\draw [dashed,double,thick] (a) -- (b) -- (c);
\draw [dashed] (a) .. controls +(30:1) and +(150:1) .. (c); 
\end{tikzpicture}$, $\begin{tikzpicture}
\node (a) at (0,0) [fill=black,draw,circle] {};
\node (b) at (1,0) [fill=white,inner sep=0pt,draw,circle] {$\pm$};
\node (c) at (2,0) [fill=black,draw,circle] {};
\draw [dashed] (a) -- (b) -- (c);
\draw [dashed,double,thick] (a) .. controls +(30:1) and +(150:1) .. (c); 
\end{tikzpicture}$ or $\begin{tikzpicture}
\node (a) at (0,0) [fill=black,draw,circle] {};
\node (b) at (1,0) [fill=black,draw,circle] {};
\node (c) at (2,0) [fill=black,draw,circle] {};
\draw [dashed,double,thick] (a) -- (b) -- (c);
\draw [dashed,double,thick] (a) .. controls +(30:1) and +(150:1) .. (c); 
\end{tikzpicture}$ . Thus, by interlacing, no cyclotomic $\mathcal{L}$-graph induces such a cycle as a subgraph.
\end{lemma}

\subsubsection{Paths with more than two consecutive weight 2 edges}

\begin{lemma}\label{maxis4cyc_E} The only cyclotomic $\mathcal{L}$ graphs of the form 
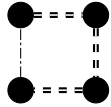
\[
\begin{tikzpicture}
\node (a) at (0,0) [fill=black,draw,circle] {};
\node (b) at (1,0) [fill=black,draw,circle] {};
\node (c) at (0,-1) [fill=black,draw,circle] {};
\node (d) at (1,-1) [fill=black,draw,circle] {};
\draw [dashed,double,thick] (a) -- (b) -- (d) -- (c);
\draw [dashdot] (c) -- (a);
\end{tikzpicture}
\] 
are equivalent to $T_4^4$ or ${T_4^4}'$ (the $k=2$ case of $T_{2k}^4$ and ${T_{2k}^4}'$ as  given in Figs. \ref{T2k4figure} and \ref{T2k4'figure}). Since in such an $\mathcal{L}$-graph all vertices have weighted degree 4, no larger $\mathcal{L}$-graph may induce a path of three consecutive weight 2 edges; by interlacing, this ensures no longer path is possible either.
\end{lemma}

\subsection{Isolated weight 2 edges}

Let $G$ be an $\mathcal{L}$-graph inducing a path $H$ with edges of weight 1, then 2, then 1. By Lemma \ref{maxis4cyc_A}, that path is of form $\begin{tikzpicture}
\node (a) at (0,0) [fill=white,inner sep=1pt,draw,shape=circle] {$*$};
\node (b) at (1,0) [fill=black,draw,circle] {};
\node (c) at (2,0) [fill=black,draw,circle] {};
\node (d) at (3,0) [fill=white,inner sep=1pt,draw,shape=circle] {$*$};
\draw [dashed] (a) -- (b) (c) -- (d);
\draw [dashed,double,thick] (b) -- (c);
\end{tikzpicture}$. However, no charged path of form $\begin{tikzpicture}
\node (a) at (0,0) [fill=white,inner sep=0pt,draw,shape=circle] {$\pm$};
\node (b) at (1,0) [fill=black,draw,circle] {};
\node (c) at (2,0) [fill=black,draw,circle] {};
\node (d) at (3,0) [fill=white,inner sep=1pt,draw,shape=circle] {$*$};
\draw [dashed] (a) -- (b) (c) -- (d);
\draw [dashed,double,thick] (b) -- (c);
\end{tikzpicture}$ is cyclotomic, so all four vertices of $H$ must be uncharged.

\begin{lemma}\label{maxis4cyc_C}
If a maximal cyclotomic $\mathcal{L}$-graph $G$ induces a subgraph of form 
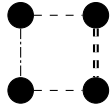
\[\begin{tikzpicture}
\node (a) at (0,0) [fill=black,draw,shape=circle] {};
\node (b) at (1,0) [fill=black,draw,circle] {};
\node (c) at (0,-1) [fill=black,draw,circle] {};
\node (d) at (1,-1) [fill=black,draw,circle] {};
\draw [dashed] (a) -- (b) (c) -- (d);
\draw [dashed,double,thick] (b) -- (d);
\draw [dashdot] (a) -- (c);
\end{tikzpicture}\] 
then $G$ is equivalent to either $S_8^*$ as given in Fig. \ref{S8*figure} or ($d=-7$ only) $S_6^\dag$ as given in Fig. \ref{S6dagfigure}.
\end{lemma}

\begin{proof} Growing terminates after four rounds and confirms that such a $G$ has either 6 or 8 vertices, and is of claimed form. In the 6 vertex case, we fix edges by switching and test the remaining possibilities for cyclotomicity; there are only two suitable choices for the remaining edge labels, $S_6^\dag$ and a graph which is confirmed to be equivalent under switching. 

In the 8 vertex case, we fix edges by switching and determine that there is only one possible set of edge labels on a 6-vertex subgraph that gives a cyclotomic subgraph (directly testing all possible combinations of unspecified labels is impractical). By interlacing, this allows us to fix those labels and test the remaining candidates; the only cyclotomic examples are equivalent to the representatives given in Fig. \ref{S8*figure}.\end{proof}

\begin{lemma}\label{maxis4cyc_B}
If a maximal cyclotomic $\mathcal{L}$-graph $G$ induces a subgraph of form $\begin{tikzpicture}
\node (a) at (0,0) [fill=black,draw,circle] {};
\node (b) at (1,0) [fill=black,draw,circle] {};
\node (c) at (2,0) [fill=black,draw,circle] {};
\draw [dashed,double,thick] (b) -- (c);
\draw [dashed] (b) -- (a) .. controls +(30:1) and +(150:1) .. (c); 
\end{tikzpicture}$ then $G$ is equivalent to the $\mathcal{L}$-graph $S_4^*$ given in Fig. \ref{S4w2figure}.
\end{lemma}
\begin{proof} No such $\mathcal{L}$-graph is cyclotomic for $d=-7$; for $d=-2$ growing terminates after a single round, and all cyclotomic examples are easily confirmed to be equivalent to $S_4^*$.\end{proof}

\section{Infinite Families of $\mathcal{L}$-graphs}\label{section6}
We have shown in the previous section that any maximal cyclotomic $\mathcal{L}$-graph neither of form $S_2, S_2', S_2^*, S_4, S_4', S_4^*, S_6^\dag, S_8^*$ nor a charged signed graph must have all edge labels from $\mathcal{L}_2\cup\mathcal{L}_1\cup\{0\}$ with at least one edge of weight 2; but any edges of weight 2 must appear in isolated pairs. The graphs $T_{2k}^4, {T_{2k}^4}'$ and  $C_{2k}^{2+}$ all satisfy these conditions; it remains to show that any maximal cyclotomic $\mathcal{L}$-graph with such properties is equivalent to one of these.

To do so, we will first demonstrate that for $d=-2,-7$ a sufficient condition for being maximally cyclotomic - that all vertices have weighted degree 4, which we describe as \emph{$4$-cyclotomic} - is also necessary. With this extra constraint, we are then able to show that any non-sporadic $\mathcal{L}$-graph is of the same form as some $T_{2k}^4, {T_{2k}^4}'$ and  $C_{2k}^{2+}$, and prove that these are representatives up to equivalence.

\subsection{Maximal cyclotomic $\mathcal{L}$-graphs are $4$-cyclotomic}
\begin{theorem}\label{neg27nonmax}
Let $G$ be a cyclotomic $\mathcal{L}$-graph with edge labels from $\mathcal{O}_{\Q(\sqrt{-2})}$ or $\mathcal{O}_{\Q(\sqrt{-7})}$. If $G$ has a vertex of weighted degree 1,2 or 3, then $G$ is nonmaximal.
\end{theorem}

\subsubsection{Excluded Subgraphs}
We identify various cyclotomic $\mathcal{L}$-graphs $H$ such that if $G$ is cyclotomic but not 4-cyclotomic and induces $H$ as a subgraph, then $G$ is not maximal. This holds when, as in the previous section, such an $H$ is (by growing) contained in only finitely many cyclotomic $\mathcal{L}$-graphs, and each of these is contained in a maximal 4-cyclotomic example; $G$ is necessarily also a proper subgraph of one of those maximal examples.

\begin{lemma}\label{nonMaxSubGraphs}
A cyclotomic $\mathcal{L}$-graph $G$ with not all vertices weight 4 is nonmaximal if it induces as subgraph either (a) an uncharged triangle or (b) a single-charged triangle:
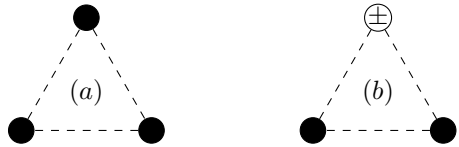
\[
\mbox{$\begin{tikzpicture}[]
\node (a) at (90:1)  [fill=black,draw,shape=circle] {};
\node (b) at (210:1)  [fill=black,draw,shape=circle] {};
\node (c) at (330:1)  [fill=black,draw,shape=circle] {};
\draw [dashed] (b) -- (a) -- (c) -- (b);
\node (name) at (0,0) {$(a)$};
\end{tikzpicture}$}
\mbox{\hspace{5em}}
\mbox{$\begin{tikzpicture}[]
\node (a) at (90:1)  [fill=white,inner sep=0pt,draw,shape=circle] {$\pm$};
\node (b) at (210:1)  [fill=black,draw,shape=circle] {};
\node (c) at (330:1)  [fill=black,draw,shape=circle] {};
\draw [dashed] (b) -- (a) -- (c) -- (b);
\node (name) at (0,0) {$(b)$};
\end{tikzpicture}$}
\]
\end{lemma}

\begin{lemma}\label{nonMaxSubGraphsw2}
A cyclotomic $\mathcal{L}$-graph $G$ with not all vertices weight 4 is nonmaximal if it induces a subgraph of any of the following forms (where cyclotomic):
\begin{itemize}
\item[(A)] Vertex with a charge and a weight 2 edge
\[
\begin{tikzpicture}[]
\node (a) at (0,0)  [fill=white,inner sep=0pt,draw,shape=circle] {$\pm$};
\node (b) at (1,0)  [fill=white,inner sep=1pt,draw,shape=circle] {$*$};
\draw [dashed,double] (b) -- (a);
\end{tikzpicture}
\]

\item[(B)] $\mathcal{L}_1,\mathcal{L}_2,\mathcal{L}_1$ Cycles
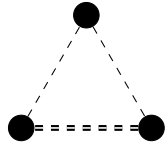
\[
\begin{tikzpicture}[]
\node (a) at (90:1)  [fill=black,draw,shape=circle] {};
\node (b) at (210:1)  [fill=black,draw,shape=circle] {};
\node (c) at (330:1)  [fill=black,draw,shape=circle] {};
\draw [dashed] (b) -- (a) -- (c);
\draw [dashed,double,thick] (c) -- (b);
\end{tikzpicture}
\]

\item[(C)] $\mathcal{L}_1,\mathcal{L}_2,\mathcal{L}_1$ Subpaths
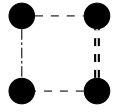
\[
\begin{tikzpicture}[]
\node (a) at (0,0)  [fill=black,draw,shape=circle] {};
\node (b) at (1,0)  [fill=black,draw,shape=circle] {};
\node (c) at (1,-1)  [fill=black,draw,shape=circle] {};
\node (d) at (0,-1)  [fill=black,draw,shape=circle] {};
\draw [dashed] (a) -- (b) (c) -- (d);
\draw [dashed,double,thick] (c) -- (b);
\draw [dashdot] (a) -- (d);
\end{tikzpicture}
\]

\item[(D)] $\mathcal{L}_2,\mathcal{L}_2,\mathcal{L}_2$ Subpaths
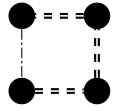
\[
\begin{tikzpicture}[]
\node (a) at (0,0)  [fill=black,draw,shape=circle] {};
\node (b) at (1,0)  [fill=black,draw,shape=circle] {};
\node (c) at (1,-1)  [fill=black,draw,shape=circle] {};
\node (d) at (0,-1)  [fill=black,draw,shape=circle] {};
\draw [dashed,double,thick] (d) -- (c) -- (b) -- (a);
\draw [dashdot] (a) -- (d);
\end{tikzpicture}
\]

\item[(E)] $\mathcal{L}_2,\mathcal{L}_1$ charged path of form
\[\begin{tikzpicture}[]
\node (x) at (2,0)  [fill=white,inner sep=0pt,draw,shape=circle] {$\pm$};
\node (a) at (1,0)  [fill=black,draw,shape=circle] {};
\node (b) at (0,0)  [fill=black,draw,shape=circle] {};
\draw [dashed] (x) -- (a);
\draw [dashed,double] (a) -- (b);
\node (a') at (1,0)  [fill=black,draw,shape=circle] {};
\end{tikzpicture}\]
\end{itemize}
\end{lemma}

\begin{proof}
\textbf (a), (b), (A) and (E) hold by growing from the given seeds, terminating with finitely many graphs each equivalent to one of the cases given in Theorems \ref{McSmThm1}, \ref{McSmThm2}, \ref{cyc_class_neg2} or \ref{cyc_class_neg7} as required; (B), (C) and (D) follow from Lemmata \ref{maxis4cyc_B}, \ref{maxis4cyc_C}, \ref{maxis4cyc_E} respectively.
\end{proof}

\subsubsection{Gram Vector Constructions}

For vectors $x,y\in \C^n$ we take as standard inner product $\ip{x}{y} = \sum_{i=1}^n x_i\overline{y_i}$ .

\begin{definition} For an $n\times n$ Hermitian matrix $A$ we describe a set $W=\{w_1,\cdots w_n\}$ as a \emph{set of Gram vectors} for $A$ if $\ip{w_i}{w_j}=A_{ij}$ for all $1\le i,j\le n$.
\end{definition}

\begin{lemma}{(Special case of \cite{Horn} Thm. 7.2.6)}\label{matrixsquareroot} Let $A$ be a positive semidefinite Hermitian matrix. Then there exists a positive semidefinite Hermitian matrix $B$ such that $B^2=A$.\end{lemma}

\begin{proposition}\label{Rhasgram} Let $A$ be an $n\times n$ positive semidefinite Hermitian $R$-matrix. Then there exists a set of Gram vectors for $A$.
\end{proposition}

Let $M$ be a matrix representative of a connected cyclotomic $\mathcal{L}$-graph $G$. Then both $A=M+2I$ and $B=(-M)+2I$ are positive semidefinite. Hence (by Proposition \ref{Rhasgram}) for a given ordering on the vertices there exist sets of Gram vectors $W$ and $W'$ for $A$ and $B$ respectively, whereby $A_{ij}=\ip{w_i}{w_j}$ and $B_{ij}=\ip{w'_i}{w'_j}$. We then have:
\begin{itemize}
\item For all $i\neq j$, $\ip{w_i}{w_j}$ and $\ip{w'_i}{w'_j}$ are in $\mathcal{L}$, with $\ip{w_i}{w_j}=-\ip{w'_i}{w'_j}$.

\item $\ip{w_i}{w_j}$ gives the label $e_{ij}$ of the edge from vertex $i$ to $j$ ($0$ if no edge); so $\ip{w_j}{w_i}=e_{ji}=\overline{e_{ij}}$ as required.
\item For all $w\in W$ and $w'\in W'$, $\ip{w}{w}$ and $\ip{w'}{w'}$ are in $\{1,2,3\}$; $\ip{w_i}{w_i}-2$ gives the charge on vertex $i$.
\item For all $i$, $\ip{w_i'}{w_i'} = 4- \ip{w_i}{w_i}$.
\end{itemize}

\begin{proposition}\label{getGramNotMax2}
Let $M$ be a matrix representative of a cyclotomic $\mathcal{L}$-graph $G$. Fix an ordered vertex labelling then determine Gram vectors $W,W'$ as above. If there exist vectors $x$, $x'$ with the following properties:
\begin{itemize}
\item $\ip{x}{x} \in \{1,2,3\}$
\item $\mbox{For all } w_i\in W$, $\ip{x}{w_i} \in \mathcal{L}$
\item There exists $w_i\in W$ such that $\ip{x}{w_i}\neq 0$
\item $\ip{x'}{x'} = 4-\ip{x}{x}$
\item $\mbox{For all } i$, $\ip{x'}{w'_i} = -\ip{x}{w_i}$
\end{itemize}
then define $A^*$ to be the matrix determined by the set of Gram vectors $W\cup\{x\}$. $M^*=A^*-2I$ is then a matrix representative of a cyclotomic $\mathcal{L}$-graph $G^*$ inducing $G$ as a proper subgraph, so $G$ is nonmaximal.
\end{proposition}

\begin{proof}
By construction $A^*$ is Hermitian and positive semidefinite. Thus $M^*$ has all eigenvalues in $[-2,\infty)$. By the first two conditions on $w$, $M^*$ has all entries in $\mathcal{L}$ so describes an $\mathcal{L}$-graph $G^*$ and by choice of Gram vectors this is an extension of $G$ by a single vertex. By the third condition $G^*$ is connected so $G$ is a proper subgraph of $G^*$; $G$ is therefore nonmaximal provided $G^*$ is cyclotomic.\\
\noindent Consider $B^*$ the Gram matrix corresponding to vectors $W'\cup\{x'\}$; by the properties of $W,W'$ and the final two conditions, $B^*$ is precisely the matrix $(-M^*)+2I$. As $B^*$ is positive semidefinite, $-M^*$ has all eigenvalues in $[-2,\infty)$. Hence $M^*$ has all eigenvalues in $(-\infty, 2]$; combined with the earlier bound this ensures all eigenvalues of $M^*$ are in $[-2,2]$ and $G^*$ is thus cyclotomic.
\end{proof}

\subsubsection{Non-maximality Proofs}

Combining the ideas of the previous two sections, we may identify cases in which a vertex of degree less than four ensures non-maximality. To complete the proof of Theorem \ref{neg27nonmax}, it is then sufficient to reduce to one of these cases.

\begin{lemma}\label{gram1} Let $G$ be a cyclotomic $\mathcal{L}$-graph containing a vertex $v$ of weight 3 such that the subgraph $H$ induced on $v$ and its neighbours is of the form 
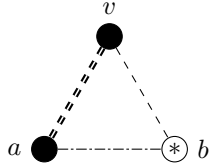
\[\begin{tikzpicture}[]
\node (x) at (90:1)  [fill=black,draw,shape=circle,label=above:{$v$}] {};
\node (a) at (210:1)  [fill=black,draw,shape=circle,label=left:{$a$}] {};
\node (b) at (330:1)  [fill=white,inner sep=1pt,draw,shape=circle,label=right:{$b$}] {$*$};
\draw [dashed] (x) -- (b);
\draw [dashed,thick,double] (x) -- (a);
\draw [dashdot] (a) -- (b);
\end{tikzpicture}\] Then $G$ is nonmaximal. 
\end{lemma}

\begin{proof}
Vertex $a$ is necessarily uncharged by Lemma \ref{nonMaxSubGraphsw2} (A). If $e_{ab}\in\mathcal{L}_1$ then (if $b$ uncharged) $G$ is nonmaximal by Lemma \ref{nonMaxSubGraphsw2} (B) or (if $b$ charged) noncyclotomic by Lemma \ref{maxis4cyc_D}. If $e_{ab}\in\mathcal{L}_2$ then $b$ is uncharged by (A), but then $G$ is noncyclotomic by Lemma \ref{maxis4cyc_D}.

Thus we conclude that $e_{ab}=0$. If $b$ is charged we have a $\mathcal{L}_2,\mathcal{L}_1$ charged path, and $G$ is nonmaximal by Lemma \ref{nonMaxSubGraphsw2} (E). Therefore $b$ is uncharged and, fixing a vertex ordering such that $v<a<b$, we have that (up to equivalence) $H$ is 

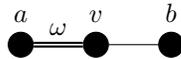
\[\begin{tikzpicture}[]
\node (x) at (1,0)  [fill=black,draw,shape=circle,label=above:{$v$}] {};
\node (a) at (0,0)  [fill=black,draw,shape=circle,label=above:{$a$}] {};
\node (b) at (2,0)  [fill=black,draw,shape=circle,label=above:{$b$}] {};
\draw (x) -- (b);
\draw [thick,double] (a) -- node[above] {$\omega$} (x);
\end{tikzpicture}\]

where $\omega=\sqrt{-2}$ or $\frac{1}{2}+\frac{\sqrt{-7}}{2}$ for $d=-2,-7$ respectively and edge labels indicate $e_{ij}$ for $i<j$ (so here $e_{va}=\omega$).

Let $W$ be the set of Gram vectors for $A=M+2I$ where $M$ is a matrix representative of $G$ with subgraph on $v,a,b$ as above. Identifying vertex $i$ with its Gram vector $w_i$, the following conditions on $W$ hold:
\[ \ip{w_v}{w_v} = \ip{w_a}{w_a} = \ip{w_b}{w_b} =2 \]
\[ \ip{w_v}{w_a} = \omega \,,\, \ip{w_v}{w_b} = 1\]
Setting $x=2w_v-\omega w_a -w_b$ we have
\[
\ip{x}{w_v} = 1,\, 
\ip{x}{w_a} = 0,\,
\ip{x}{w_b} = 0,\,
\ip{x}{x} = 2
\]
Further, for any $w_i\in W\backslash\{w_v,w_a,w_b\}$ $\ip{w_v}{w_i}=0$ by assumption so
\[\ip{x}{w_i}=-\omega\ip{w_a}{w_i}-\ip{w_b}{w_i}\]
but (fixing $v<a<b<i$) testing confirms that the subgraph induced on $v,a,b,i$
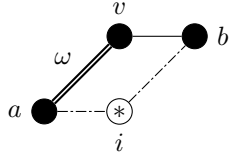
\[\begin{tikzpicture}[]
\node (x) at (0,1)  [fill=black,draw,shape=circle,label=above:{$v$}] {};
\node (a) at (-1,0)  [fill=black,draw,shape=circle,label=left:{$a$}] {};
\node (b) at (1,1)  [fill=black,draw,shape=circle,label=right:{$b$}] {};
\node (i) at (0,0) [fill=white,inner sep=1pt,draw,shape=circle,label=below:{$i$}] {$*$};
\draw (x) -- (b);
\draw [thick,double] (a) -- node[above left] {$\omega$} (x);
\draw [dashdot] (a) -- (i) -- (b);
\end{tikzpicture}\]
is cyclotomic only if $-\omega\ip{w_a}{w_i}-\ip{w_b}{w_i}\in\mathcal{L}$; thus $\ip{x}{w_i}\in\mathcal{L}$ $\mbox{for all } w_i\in W$, $\ip{x}{x}\in\{1,2,3\}$ and $\ip{x}{w_v}\neq 0$. So all conditions on $x$ required by Proposition  \ref{getGramNotMax2} hold.

With the same vertex labelling and ordering we now consider $W'$ the Gram vectors of $B=(-M)+2I$, for which we have the following:
\[ \ip{w_v'}{w_v'} = \ip{w_a'}{w_a'} = \ip{w_b'}{w_b'} =2 \]
\[ \ip{w_v'}{w_a'} = -\omega \,,\, \ip{w_v'}{w_b'} = -1\]
Setting $x'=-2w_v'-\omega w_a'-w_b'$ we have $\ip{x'}{w_i'}=-\ip{x}{w_i}$ for $w_i'\in\{w_v',w_a',w_b'\}$, $\ip{x'}{x'}=4-\ip{x}{x}$ and \[\ip{x}{w_i'}=-\omega\ip{w_a'}{w_i'}-\ip{w_b'}{w_i'} = \omega\ip{w_a}{w_i} + \ip{w_b}{w_i} = -\ip{x}{w_i}\]
for $w_i'\in W'\backslash\{w_v',w_a',w_b'\}$, so by Proposition  \ref{getGramNotMax2} $G$ is nonmaximal.
\end{proof}

The remaining cases are similar, so we omit some of the details; full versions can be found in \cite{Taylor}.

\begin{lemma}\label{gram2}
Let $G$ be a cyclotomic $\mathcal{L}$-graph containing a vertex $v$ of weight 3 such that the subgraph $H$ induced on $v$ and its neighbours is of the form
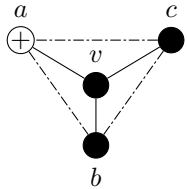
\[\begin{tikzpicture}[]
\node (c) at (2,0)  [fill=black,draw,shape=circle,label=above:{$c$}] {}; 
\node (x) at (1,-0.6)  [fill=black,draw,shape=circle,label=above:{$v$}] {}; 
\node (a) at (0,0)  [fill=white,inner sep=0pt,draw,shape=circle,label=above:{$a$}] {$+$}; 
\node (b) at (1,-1.4) [fill=black,draw,shape=circle,label=below:{$b$}] {};
\draw (a) -- (x) -- (b) (x) -- (c);
\draw [dashdot] (a) -- (b) -- (c) -- (a);
\end{tikzpicture}
\]
Then $G$ is nonmaximal.
\end{lemma}

\begin{proof}
By Lemma \ref{nonMaxSubGraphsw2} (A) we have $e_{ab},e_{ac}\not\in\mathcal{L}_2$; $e_{bc}\not\in\mathcal{L}_2$ by part (B) of the same. Further, $e_{ab},e_{ac}\not\in\mathcal{L}_1$ by Lemma \ref{nonMaxSubGraphs} (b) and $e_{bc}\not\in\mathcal{L}_1$ by part (a) of the same. So $e_{ab}=e_{ac}=e_{bc}=0$ and thus we have that $H$ is 
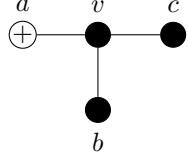
\[\begin{tikzpicture}[]
\node (c) at (2,0)  [fill=black,draw,shape=circle,label=above:{$c$}] {}; 
\node (x) at (1,0)  [fill=black,draw,shape=circle,label=above:{$v$}] {}; 
\node (a) at (0,0)  [fill=white,inner sep=0pt,draw,shape=circle,label=above:{$a$}] {$+$}; 
\node (b) at (1,-1) [fill=black,draw,shape=circle,label=below:{$b$}] {};
\draw (a) -- (x) -- (b) (x) -- (c);
\end{tikzpicture}\]

Let $W$ be the set of Gram vectors for $A=M+2I$ where $M$ is a matrix representative of $G$ with the subgraph on $v,a,b,c$ as above. Identifying vertex $i$ with its Gram vector $w_i$ and setting $x=2w_v-w_a-w_b-w_c$ we then have
\[
\ip{x}{w_a} = -1 \\ 
\ip{x}{w_b} = 0 \,,
\ip{x}{w_c} = 0 \,,
\ip{x}{w_v} = 1 \,,
\ip{x}{x} = 3 
\]
Further, if $w_i\in W\backslash \{w_v,w_a,w_b,w_c\}$ then by assumption $\ip{w_v}{w_i}=0$, so \[\ip{x}{w_i} = -(\ip{w_a}{w_i} + \ip{w_b}{w_i} + \ip{w_c}{w_i})\] and the subgraph induced on $v,a,b,c,i$ is of form  
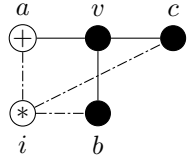
\[\begin{tikzpicture}[]
\node (c) at (2,0)  [fill=black,draw,shape=circle,label=above:{$c$}] {}; 
\node (x) at (1,0)  [fill=black,draw,shape=circle,label=above:{$v$}] {}; 
\node (a) at (0,0)  [fill=white,inner sep=0pt,draw,shape=circle,label=above:{$a$}] {$+$}; 
\node (b) at (1,-1) [fill=black,draw,shape=circle,label=below:{$b$}] {};
\node (i) at (0,-1) [fill=white,inner sep=1pt,draw,shape=circle,label=below:{$i$}] {$*$}; 
\draw (a) -- (x) -- (b) (x) -- (c);
\draw [dashdot] (a) -- (i) -- (b) (i) -- (c);
\end{tikzpicture}\] which (fixing an ordering) is cyclotomic only if $\ip{w_a}{w_i}+\ip{w_b}{w_i}+\ip{w_c}{w_i}\in\mathcal{L}$. 

Thus $\ip{x}{x}\in\{1,2,3\}$ and $\ip{x}{w_i}\in\mathcal{L}$  $\mbox{for all } w_i\in W$. Further, $\ip{x}{w_v} \neq 0$. With the same vertex labelling we now consider $W'$ the Gram vectors of $B=(-M)+2I$; setting $x'=-w_a'-w_b'-w_c'-2w_v'$ we then have $\ip{x'}{w_i'}=-\ip{x}{w_i}$ for $w_i'\in\{w_a',w_b',w_c',w_v'\}$, $\ip{x'}{x'}=4-\ip{x}{x}$. 
Further, if $w_i'\in W'\backslash \{w_v',w_a',w_b',w_c'\}$ then by assumption $\ip{w_v'}{w_i'}=0$ and so \[\ip{x'}{w_i'} =  -\ip{w_a'}{w_i'} - \ip{w_b'}{w_i'} -\ip{w_c'}{w_i'} = \ip{w_a}{w_i} + \ip{w_b}{w_i} + \ip{w_c}{w_i} = -\ip{x}{w_i}\]
Thus $G$ is nonmaximal by Proposition  \ref{getGramNotMax2}.
\end{proof}

\begin{lemma}\label{gram3}
Let $G$ be a cyclotomic $\mathcal{L}$-graph containing a vertex $v$ of weight 3 such that the subgraph $H$ induced on $v$ and its neighbours is of the form
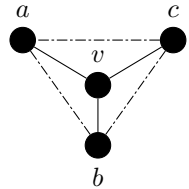
\[\begin{tikzpicture}[]
\node (c) at (2,0)  [fill=black,draw,shape=circle,label=above:{$c$}] {}; 
\node (x) at (1,-0.6)  [fill=black,draw,shape=circle,label=above:{$v$}] {}; 
\node (a) at (0,0)  [fill=black,draw,shape=circle,label=above:{$a$}] {}; 
\node (b) at (1,-1.4) [fill=black,draw,shape=circle,label=below:{$b$}] {};
\draw (a) -- (x) -- (b) (x) -- (c);
\draw [dashdot] (a) -- (b) -- (c) -- (a);
\end{tikzpicture}
\]
Then $G$ is nonmaximal.
\end{lemma}

\begin{proof}
By Lemma \ref{nonMaxSubGraphsw2} (B), $e_{ab},e_{ac},e_{bc}\not\in\mathcal{L}_2$ ; further, $e_{ab},e_{ac},e_{bc}\not\in\mathcal{L}_1$ by Lemma \ref{nonMaxSubGraphs} (a). So $e_{ab}=e_{ac}=e_{bc}=0$ and thus we have that $H$ is 
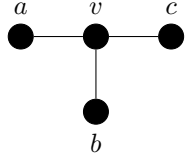
\[\begin{tikzpicture}[]
\node (c) at (2,0)  [fill=black,draw,shape=circle,label=above:{$c$}] {}; 
\node (x) at (1,0)  [fill=black,draw,shape=circle,label=above:{$v$}] {}; 
\node (a) at (0,0)  [fill=black,draw,shape=circle,label=above:{$a$}] {}; 
\node (b) at (1,-1) [fill=black,draw,shape=circle,label=below:{$b$}] {};
\draw (a) -- (x) -- (b) (x) -- (c);
\end{tikzpicture}\]

Let $W$ be the set of Gram vectors for $A=M+2I$ where $M$ is a matrix representative of $G$ with the subgraph on $v,a,b,c$ as above. Identifying vertex $i$ with its Gram vector $w_i$and setting $x=2w_v-w_a-w_b-w_c$ we then have
\[
\ip{x}{w_a} = \ip{x}{w_b} = \ip{x}{w_c} = 0\,,
\ip{x}{w_v} = 1\,, \ip{x}{x} = 2 
\]
Further, if $w_i\in W\backslash \{w_v,w_a,w_b,w_c\}$ then by assumption $\ip{w_v}{w_i}=0$ and so \[\ip{x}{w_i} = -(\ip{w_a}{w_i} + \ip{w_b}{w_i} + \ip{w_c}{w_i})\] and the subgraph induced on $v,a,b,c,i$ is of form
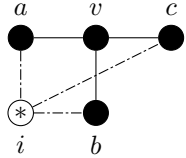
\[\begin{tikzpicture}[]
\node (c) at (2,0)  [fill=black,draw,shape=circle,label=above:{$c$}] {}; 
\node (x) at (1,0)  [fill=black,draw,shape=circle,label=above:{$v$}] {}; 
\node (a) at (0,0)  [fill=black,draw,shape=circle,label=above:{$a$}] {}; 
\node (b) at (1,-1) [fill=black,draw,shape=circle,label=below:{$b$}] {};
\node (i) at (0,-1) [fill=white,inner sep=1pt,draw,shape=circle,label=below:{$i$}] {$*$}; 
\draw (a) -- (x) -- (b) (x) -- (c);
\draw [dashdot] (a) -- (i) -- (b) (i) -- (c);
\end{tikzpicture}\] which (fixing an ordering) is cyclotomic only if $\ip{w_a}{w_i}+\ip{w_b}{w_i}+\ip{w_c}{w_i}\in\mathcal{L}$. 

Thus $x\cdot x\in\{1,2,3\}$ and $x\cdot w_i\in\mathcal{L}$  $\mbox{for all } w_i\in W$. Further, $\ip{x}{w_v} \neq 0$ With the same vertex labelling we now consider $W'$ the Gram vectors of $B=(-M)+2I$; setting $x'=-w_a'-w_b'-w_c'-2w_v'$ we then have $\ip{x'}{w_i'}=-\ip{x}{w_i}$ for $w_i'\in\{w_v',w_a',w_b',w_c'\}$, and $\ip{x'}{x'}=4-\ip{x}{x}$. Further, if $w_i'\in W'\backslash \{w_v',w_a',w_b',w_c'\}$ then by assumption $\ip{w_v'}{w_i'}=0$,thus \[\ip{x'}{w_i'} =  -\ip{w_a'}{w_i'} - \ip{w_b'}{w_i'} -\ip{w_c'}{w_i'} = \ip{w_a}{w_i} + \ip{w_b}{w_i} + \ip{w_c}{w_i} = -\ip{x}{w_i}\]
Thus $G$ is nonmaximal by Proposition  \ref{getGramNotMax2}.
\end{proof}

\begin{lemma}\label{gram4} 
Let $G$ be a cyclotomic $\mathcal{L}$-graph containing a vertex $v$ of weight 2 in a subgraph $H$ of the form
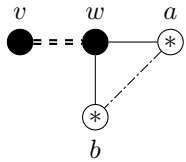
\[\begin{tikzpicture}[]
\node (w) at (1,0)  [fill=black,draw,shape=circle,label=above:{$w$}] {}; 
\node (v) at (0,0) [fill=black,draw,shape=circle,label=above:{$v$}] {};
\node (a) at (2,0)  [fill=white,inner sep=1pt,draw,shape=circle,label=above:{$a$}] {$*$};
\node (b) at (1,-1) [fill=white,inner sep=1pt,draw,shape=circle,label=below:{$b$}] {$*$};
\draw [dashed,double,thick] (v) -- (w);
\draw (a) -- (w) -- (b);
\draw [dashdot] (a) -- (b);
\end{tikzpicture}
\]
Then $G$ is nonmaximal.
\end{lemma}

\begin{proof}
 By Lemma \ref{nonMaxSubGraphsw2} (E), $G$ is nonmaximal if either $a$ or $b$ is charged. But if neither is charged then $e_{ab}\not\in\mathcal{L}_2$ by part (B) of the same and $e_{ab}\not\in\mathcal{L}_1$ by Lemma \ref{nonMaxSubGraphs} (a). Hence we may assume $e_{ab}=0$ and that $a,b$ are uncharged; fixing an ordering $v<w<a<b$ we have that $H$ is, up to equivalence, 
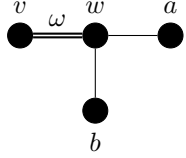
\[
\begin{tikzpicture}[]
\node (w) at (1,0)  [fill=black,draw,shape=circle,label=above:{$w$}] {}; 
\node (v) at (0,0) [fill=black,draw,shape=circle,label=above:{$v$}] {};
\node (a) at (2,0)  [fill=black,draw,shape=circle,label=above:{$a$}] {};
\node (b) at (1,-1) [fill=black,draw,shape=circle,label=below:{$b$}] {};
\draw [thick,double] (v) -- node[above] {$\omega$} (w);
\draw (a) -- (w) -- (b);
\end{tikzpicture}
\]
where $\omega=\sqrt{-2}$ or $\frac{1}{2}+\frac{\sqrt{-7}}{2}$ for $d=-2,-7$ respectively.

Let $W$ be the set of Gram vectors for $A=M+2I$ where $M$ is a matrix representative of $G$ with subgraph on $v,w,a,b$ as above. Identifying vertex $i$ with its Gram vector $w_i$ and setting $x=w_w-w_a-w_b$ we have
\[
\ip{x}{w_v} = \overline{\omega}\,,
\ip{x}{w_w} = 0\,,
\ip{x}{w_a} = -1\,,
\ip{x}{w_b} = -1\,,
\ip{x}{x} = 2\]
Further, for any $w_i\in W\backslash\{w_v,w_w,w_a,w_b\}$, $\ip{w_w}{w_i}=0$ since $w$ has weighted degree 4, so
\[\ip{x}{w_i}=-\ip{w_a}{w_i}-\ip{w_b}{w_i}\]
but (fixing an ordering) testing confirms that for any vertex $i$, the subgraph induced on $v,w,a,b,i$
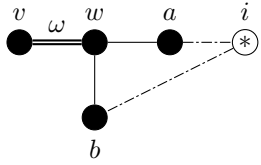
\[
\begin{tikzpicture}[]
\node (w) at (1,0)  [fill=black,draw,shape=circle,label=above:{$w$}] {}; 
\node (v) at (0,0) [fill=black,draw,shape=circle,label=above:{$v$}] {};
\node (a) at (2,0)  [fill=black,draw,shape=circle,label=above:{$a$}] {};
\node (b) at (1,-1) [fill=black,draw,shape=circle,label=below:{$b$}] {};
\node (i) at (3,0) [fill=white,inner sep=1pt,draw,shape=circle,label=above:{$i$}] {$*$};
\draw [thick,double] (v) -- node[above] {$\omega$} (w);
\draw (a) -- (w) -- (b);
\draw [dashdot] (a) -- (i) -- (b);
\end{tikzpicture}
\]
is cyclotomic only if $\ip{w_a}{w_i}=-\ip{w_b}{w_i}$ and so $\ip{x}{w_i}=0$ for all such $w_i$. With the same vertex labelling and ordering we now consider $W'$ the Gram vectors of $B=(-M)+2I$; setting $x'=w_w'+w_a'+w_b'$ we have $\ip{x'}{w_i'}=-\ip{x}{w_i}$ for $w_i'\in\{w_v',w_w',w_a',w_b'\}$, $\ip{x'}{x'}=4-\ip{x}{x}$ and, for $w_i'\in W'\backslash\{w_v,w_w,w_a,w_b\}$, \[ \ip{x'}{w_i'}=\ip{w_a'}{w_i'}+\ip{w_b'}{w_i'}=-(\ip{w_a}{w_i}+\ip{w_b}{w_i})=-(0)=-\ip{x}{w_i}\] 
Thus $G$ is nonmaximal by Proposition  \ref{getGramNotMax2}.
\end{proof}

\subsubsection{Proof of Theorem \ref{neg27nonmax}}
Let $G$ be a cyclotomic $\mathcal{L}$-graph with a vertex of degree less than four, and at least one edge of weight 2. We will show that $G$ is nonmaximal. We first note the following two results, which hold by direct testing:

\begin{lemma}\label{34sub} There are no cyclotomic $\mathcal{L}$-graphs of the form
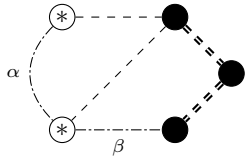
\[
\begin{tikzpicture}[scale=1.5]
\node (a0) at (0,0) [fill=black,draw,circle] {};
\node (z) at (0.5,-0.5) [fill=black,draw,circle] {};
\node (b0) at (0,-1) [fill=black,draw,circle] {};
\node (x1) at (-1,0) [fill=white,inner sep=1pt,draw,circle] {$*$};
\node (x2) at (-1,-1) [fill=white,inner sep=1pt,draw,circle] {$*$};
\draw [double,dashed,thick] (a0) -- (z) -- (b0);
\draw [dashed] (x1) -- (a0) -- (x2);
\draw [dashdot] (x2) -- node [below,black] {\scriptsize$\beta$} (b0); 
\draw [dashdot] (x1) .. controls +(225:0.5) and +(135:0.5) .. node [left,black] {\scriptsize$\alpha$} (x2); 
\end{tikzpicture}
\]
for $\alpha\in\mathcal{L},\beta\in\mathcal{L}_1\cup\{0\}$.

Thus no cyclotomic $\mathcal{L}$-graph has such a graph as an induced subgraph.
\end{lemma}

\begin{lemma}\label{k1nocyc}
 There are no cyclotomic graphs of the form
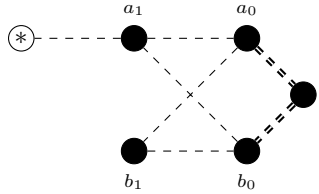
\[
\begin{tikzpicture}[scale=1.5]
\node (a0) at (0,0) [fill=black,draw,circle,label=above:{\scriptsize$a_0$}] {};
\node (z) at (0.5,-0.5) [fill=black,draw,circle] {};
\node (b0) at (0,-1) [fill=black,draw,circle,label=below:{\scriptsize$b_0$}] {};
\node (x1) at (-1,0) [fill=black,draw,circle,label=above:{\scriptsize$a_1$}] {};
\node (x2) at (-1,-1) [fill=black,draw,circle,label=below:{\scriptsize$b_1$}] {};
\node (x) at (-2,0) [fill=white,inner sep=1pt,draw,circle] {$*$};
\draw [double,dashed,thick] (a0) -- (z) -- (b0);
\draw [dashed] (x1) -- (a0) -- (x2) -- (b0) -- (x1) -- (x);
\end{tikzpicture}
\]
\end{lemma}

and the following Corollary of Theorems \ref{McSmThm1}, \ref{McSmThm2}:

\begin{corollary}\label{kgt1nocyc}
 There are no cyclotomic $\mathcal{L}$-graphs of the form
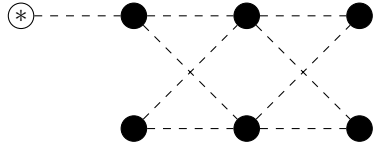
\[
\begin{tikzpicture}[scale=1.5]
\node (1) at (0,0) [fill=black,draw,circle] {};
\node (2) at (0,-1) [fill=black,draw,circle] {};
\node (3) at (1,0) [fill=black,draw,circle] {};
\node (4) at (1,-1) [fill=black,draw,circle] {};
\node (5) at (2,0) [fill=black,draw,circle] {};
\node (6) at (2,-1) [fill=black,draw,circle] {};
\node (x) at (-1,0) [fill=white,inner sep=1pt,draw,circle] {$*$};
\draw [dashed] (x) -- (1) -- (3) -- (5) (2) -- (4) -- (6) (1) -- (4) -- (5) (2) -- (3) -- (6);
\end{tikzpicture}
\]
 \end{corollary}

Let $x,y$ be vertices of $G$ joined by a weight 2 edge; by Lemma \ref{maxis4cyc_D} we may assume $x,y$ are uncharged. If $x$ and $y$ have no further neighbours then they are the entirety of $G$ which is trivially nonmaximal. 

If there are no additional weight 2 edges incident at either $x$ or $y$ then there must be a weight 1 edge incident at one but - by Lemma \ref{nonMaxSubGraphsw2} (B) and (C) - not the other. W.l.o.g, let $x$ have no other neighbours in $G$. If $y$ only has one more neighbour, then it is a weight 3 vertex satisfying the conditions of Lemma \ref{gram1}. If $y$ has two neighbours, then $x$ is a weight 2 vertex satisfying the conditions of Lemma \ref{gram4}. In either case, this ensures that $G$ is nonmaximal.

So there must be an additional edge of weight 2 incident at either $x$ or $y$ but - by Lemma \ref{maxis4cyc_D} and Lemma \ref{nonMaxSubGraphsw2} (D) - not both. W.l.o.g. let it be $x$; Lemma \ref{nonMaxSubGraphsw2} (A) ensures the vertex $z$ joined in this way to $x$ is uncharged; Lemma \ref{maxis4cyc_D} forces $e_{yz}=0$. Again, if neither $y$ nor $z$ has further neighbours in $G$, then $x,y,z$ are all the vertices of $G$ which is clearly nonmaxmimal. Further, if either $y$ or $z$ has only one neighbour, then it is a weight 3 vertex satisfying the conditions of Lemma \ref{gram1}, so $G$ is nonmaximal.

Therefore both $y$ and $z$ have two neighbours in $G$; Lemma \ref{34sub} ensures that they have common neighbours $a_1,b_1$; if either of these is charged, we necessarily have a 4-cyclotomic graph of form $C_4^{2\pm}$, which is a contradiction. Thus $a_1,b_1$ are uncharged; further, they cannot be neighbours without violating cyclotomicity, so $G$ induces a subgraph of form:
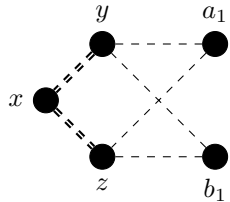
\[
\begin{tikzpicture}[scale=1.5]
\node (y) at (0,0) [fill=black,draw,circle,label=above:{$y$}] {};
\node (x) at (-0.5,-0.5) [fill=black,draw,circle,label=left:{$x$}] {};
\node (z) at (0,-1) [fill=black,draw,circle,label=below:{$z$}] {};
\node (a1) at (1,0) [fill=black,draw,circle,label=above:{$a_1$}] {};
\node (b1) at (1,-1) [fill=black,draw,circle,label=below:{$b_1$}] {};
\draw [double,dashed,thick] (y) -- (x) -- (z);
\draw [dashed] (y) -- (a1) -- (z) -- (b1) -- (y);
\end{tikzpicture}
\]

By Lemma \ref{k1nocyc} any neighbour of $a_1$ is a neighbour of $b_1$. So either we have that $G$ is the graph above (nonmaximal by embedding in $T_6^4$); that there is a single common neighbour $a_2$ of $a_1,b_1$:

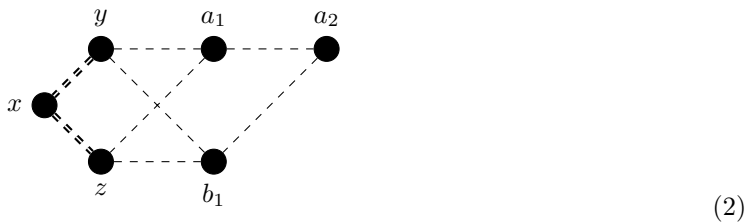
\begin{equation}\label{gram5}
\begin{tikzpicture}[scale=1.5]
\node (y) at (0,0) [fill=black,draw,circle,label=above:{$y$}] {};
\node (x) at (-0.5,-0.5) [fill=black,draw,circle,label=left:{$x$}] {};
\node (z) at (0,-1) [fill=black,draw,circle,label=below:{$z$}] {};
\node (a1) at (1,0) [fill=black,draw,circle,label=above:{$a_1$}] {};
\node (a2) at (2,0) [fill=black,draw,circle,label=above:{$a_2$}] {};
\node (b1) at (1,-1) [fill=black,draw,circle,label=below:{$b_1$}] {};
\draw [double,dashed,thick] (y) -- (x) -- (z);
\draw [dashed] (y) -- (a1) -- (z) -- (b1) -- (y) (a1) -- (a2) -- (b1);
\end{tikzpicture}
\end{equation}

or there is a pair of common neighbours $a_2,b_2$. Lemma \ref{kgt1nocyc} ensures that for each pair $a_i,b_i$ any neighbour of one is a neighbour of the other. Thus we continue to identify pairs of common neighbours until we reach a $j$ such that $a_j,b_j$ have weight less than four; (\ref{gram5}) is the case $j=1$. 

If $a_j,b_j$ are both of weight 2 then we have that $G$ is a chain of length $j$, which is nonmaximal by embedding into, for instance, a $T_{2k}^4$. Otherwise they are of weight 3; if their mutual neighbour is uncharged then $a_j$ satisfies the conditions of Lemma \ref{gram3}, but if it is charged then Lemma \ref{gram2} applies. Thus $G$ is nonmaximal; this completes the proof of Theorem \ref{neg27nonmax}.  

\section{Equivalence Classes of Infinite Families}

\begin{definition}\label{kcyclinder}
We describe any $2m$-vertex graph of the form
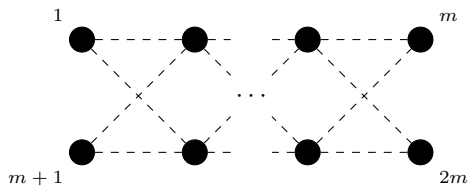
\[
\begin{tikzpicture}[scale=1.5]
\node (2) at (1,0) [fill=black,draw,circle,label=above left:{\scriptsize$1$}] {};
\node (3) at (2,0) [fill=black,draw,circle] {};
\node (3a) at (2.4,0) {};
\node (3b) at (2.4,-0.4) {};
\draw [dashed] (2) -- (3) -- (3a); 
\draw [dashed] (3) -- (3b);
\node (k+2) at (1,-1) [fill=black,draw,circle,label=below left:{\scriptsize$m+1$}] {};
\node (k+3) at (2,-1) [fill=black,draw,circle] {};
\node (k+3a) at (2.4,-1) {};
\node (k+3b) at (2.4,-0.6) {};
\draw [dashed] (k+2) -- (k+3) (2) -- (k+3) (k+2) -- (3) (k+3a) -- (k+3) -- (k+3b);
\node (spacer) at (2.5,-0.5) {$\cdots$};
\node (k-1b) at (2.6,-0.4) {};
\node (2k-1b) at (2.6,-0.6) {};
\node (k-1a) at (2.6,0) {};
\node (k-1) at (3,0) [fill=black,draw,circle] {};
\node (k) at (4,0) [fill=black,draw,circle,label=above right:{\scriptsize$m$}] {};
\draw [dashed] (k-1a) -- (k-1) -- (k);
\node (2k-1a) at (2.6,-1) {};
\node (2k-1) at (3,-1) [fill=black,draw,circle] {};
\node (2k) at (4,-1) [fill=black,draw,circle,label=below right:{\scriptsize$2m$}] {};
\draw [dashed] (2k) -- (2k-1) -- (2k-1a);
\draw [dashed] (2k-1) -- (2k-1b);
\draw [dashed] (k-1b) -- (k-1) -- (2k);
\draw [dashed] (2k-1) -- (k);
\end{tikzpicture}
\] 
as a \emph{cylinder of length $m$}.
\end{definition}

We note the following consequences of Theorems \ref{McSmThm1}, \ref{McSmThm2}:

\begin{corollary}\label{cylinderfixall}
If $g$ is a  cyclotomic $2m$-vertex cylinder of length $m\ge4$ with all edge labels $\pm1$ then $g$ is equivalent to the signed graph
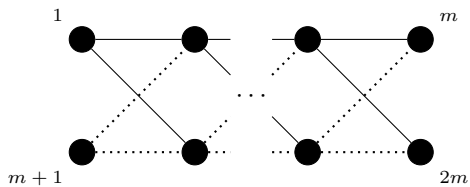
\[
\begin{tikzpicture}[scale=1.5]
\node (2) at (1,0) [fill=black,draw,circle,label=above left:{\scriptsize$1$}] {};
\node (3) at (2,0) [fill=black,draw,circle] {};
\node (3a) at (2.4,0) {};
\node (3b) at (2.4,-0.4) {};
\draw (2) -- (3) -- (3a); 
\draw (3) -- (3b);
\node (k+2) at (1,-1) [fill=black,draw,circle,label=below left:{\scriptsize$m+1$}] {};
\node (k+3) at (2,-1) [fill=black,draw,circle] {};
\node (k+3a) at (2.4,-1) {};
\node (k+3b) at (2.4,-0.6) {};
\draw [dotted,thick] (k+2) -- (k+3) -- (k+3a);
\draw [dotted,thick] (k+3) -- (k+3b);
\draw [dotted,thick] (k+2) -- (3);
\draw (2) -- (k+3);
\node (spacer) at (2.5,-0.5) {$\cdots$};
\node (k-1b) at (2.6,-0.4) {};
\node (2k-1b) at (2.6,-0.6) {};
\node (k-1a) at (2.6,0) {};
\node (k-1) at (3,0) [fill=black,draw,circle] {};
\node (k) at (4,0) [fill=black,draw,circle,label=above right:{\scriptsize$m$}] {};
\draw (k-1a) -- (k-1) -- (k);
\node (2k-1a) at (2.6,-1) {};
\node (2k-1) at (3,-1) [fill=black,draw,circle] {};
\node (2k) at (4,-1) [fill=black,draw,circle,label=below right:{\scriptsize$2m$}] {};
\draw [dotted,thick] (2k-1a) -- (2k-1) -- (2k);
\draw (2k-1b) -- (2k-1) (k-1) -- (2k);
\draw [dotted,thick] (k-1b) -- (k-1) (2k-1) -- (k);
\end{tikzpicture}
\]
\end{corollary}

\begin{lemma}\label{samesign}
A charged signed graph of the form
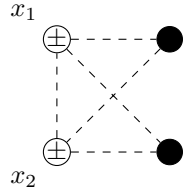
\[
\begin{tikzpicture}[scale=1.5]
\node (a) at (0,0) [fill=white,inner sep=0pt,draw,circle,label=above left:{$x_1$}] {$\pm$};
\node (b) at (1,0) [fill=black,draw,circle] {};
\node (c) at (0,-1) [fill=white,inner sep=0pt,draw,circle,label=below left:{$x_2$}] {$\pm$};
\node (d) at (1,-1) [fill=black,draw,circle] {};
\draw [dashed] (b) -- (c) -- (a) -- (d);
\draw [dashed] (a) -- (b);
\draw [dashed] (c) -- (d);
\end{tikzpicture}
\] is cyclotomic only if $x_1=x_2$ (that is, the two charged vertices have the same charge).
\end{lemma}

\begin{corollary}\label{C2k++cor}
A charged signed graph of the form 
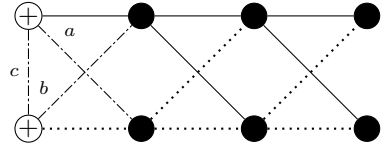
\[
\begin{tikzpicture}[scale=1.5]
\node (a) at (0,0) [fill=white,inner sep=0pt,draw,circle] {$+$};
\node (b) at (1,0) [fill=black,draw,circle] {};
\node (c) at (2,0) [fill=black,draw,circle] {};
\node (c2) at (3,0) [fill=black,draw,circle] {};
\node (d) at (0,-1) [fill=white,inner sep=0pt,draw,circle] {$+$};
\node (e) at (1,-1) [fill=black,draw,circle] {};
\node (f) at (2,-1) [fill=black,draw,circle] {};
\node (f2) at (3,-1) [fill=black,draw,circle] {};
\draw (a) -- (b) -- (c) -- (c2) (b) -- (f) (c) -- (f2);
\draw [dotted,thick] (d) -- (e) -- (f) -- (f2) (e) -- (c) (f) -- (c2);
\draw [dashdot] (a) -- node [left,black] {\scriptsize$c$} (d) (a) -- node [above left,black,yshift=1em] {\scriptsize$a$} (e) (d) -- node [below left,black,xshift=-1em] {\scriptsize$b$} (b);
\end{tikzpicture}
\] is cyclotomic if and only if $a=c=1,b=-1$, whilst a charged signed graph of the form  
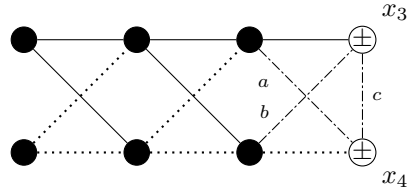
\[
\begin{tikzpicture}[scale=1.5]
\node (a) at (4,0) [fill=white,inner sep=0pt,draw,circle,label=above right:{$x_3$}] {$\pm$};
\node (b) at (1,0) [fill=black,draw,circle] {};
\node (c) at (2,0) [fill=black,draw,circle] {};
\node (c2) at (3,0) [fill=black,draw,circle] {};
\node (d) at (4,-1) [fill=white,inner sep=0pt,draw,circle,label=below right:{$x_4$}] {$\pm$};
\node (e) at (1,-1) [fill=black,draw,circle] {};
\node (f) at (2,-1) [fill=black,draw,circle] {};
\node (f2) at (3,-1) [fill=black,draw,circle] {};
\draw (b) -- (c) -- (c2) -- (a) (b) -- (f) (c) -- (f2);
\draw [dotted,thick] (e) -- (f) -- (f2) -- (d) (e) -- (c) (f) -- (c2);
\draw [dashdot] (a) -- node [right,black] {\scriptsize$c$} (d) (f2) -- node [above left,black,xshift=-1em] {\scriptsize$a$} (a) (c2) -- node [below left,black,xshift=-1em] {\scriptsize$b$} (d);
\end{tikzpicture}
\] is cyclotomic if and only if $x_3=x_4=a=1,b=c=-1$ or $x_3=x_4=b=-1,a=c=1$.
\end{corollary}

\subsection{$\mathcal{L}$-graphs of the form $\mathcal{T}^4_{2k}$}

\begin{definition}\label{T2k4def}
For $k=L+1\ge2$ define the $2k$-vertex form $\mathcal{T}_{2k}^{4}$ by
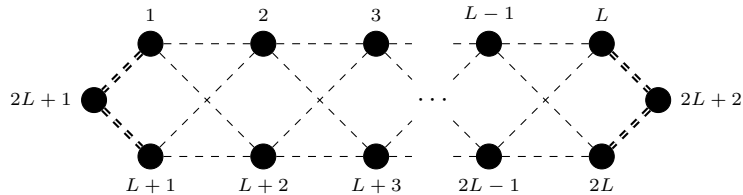
\[
\begin{tikzpicture}[scale=1.5]
\node (1) at (0,0) [fill=black,draw,circle,label=above:{\scriptsize$1$}] {};
\node (2) at (1,0) [fill=black,draw,circle,label=above:{\scriptsize$2$}] {};
\node (3) at (2,0) [fill=black,draw,circle,label=above:{\scriptsize$3$}] {};
\node (3a) at (2.4,0) {};
\node (3b) at (2.4,-0.4) {};
\draw [dashed] (1) -- (2) -- (3) -- (3a); 
\draw [dashed] (3) -- (3b);
\node (L+1) at (0,-1) [fill=black,draw,circle,label=below:{\scriptsize$L+1$}] {};
\node (L+2) at (1,-1) [fill=black,draw,circle,label=below:{\scriptsize$L+2$}] {};
\node (L+3) at (2,-1) [fill=black,draw,circle,label=below:{\scriptsize$L+3$}] {};
\node (k+3a) at (2.4,-1) {};
\node (k+3b) at (2.4,-0.6) {};
\draw [dashed] (L+1) -- (L+2) -- (L+3) -- (k+3a);
\draw [dashed] (L+3) -- (k+3b);
\draw [dashed] (1) -- (L+2) -- (3);
\draw [dashed] (L+1) -- (2) -- (L+3);
\node (spacer) at (2.5,-0.5) {$\cdots$};
\node (k-1b) at (2.6,-0.4) {};
\node (2k-1b) at (2.6,-0.6) {};
\node (k-1a) at (2.6,0) {};
\node (L-1) at (3,0) [fill=black,draw,circle,label=above:{\scriptsize$L-1$}] {};
\node (L) at (4,0) [fill=black,draw,circle,label=above:{\scriptsize$L$}] {};
\draw [dashed] (k-1a) -- (L-1) -- (L);
\node (2k-1a) at (2.6,-1) {};
\node (2L-1) at (3,-1) [fill=black,draw,circle,label=below:{\scriptsize$2L-1$}] {};
\node (2L) at (4,-1) [fill=black,draw,circle,label=below:{\scriptsize$2L$}] {};
\draw [dashed] (2k-1a) -- (2L-1) -- (2L);
\draw [dashed] (2k-1b) -- (2L-1) -- (L);
\draw [dashed] (k-1b) -- (L-1) -- (2L);
\node (2L+1) at (-0.5,-0.5) [fill=black,draw,circle,label=left:{\scriptsize$2L+1$}] {};
\node (2L+2) at (4.5,-0.5) [fill=black,draw,circle,label=right:{\scriptsize$2L+2$}] {};
\draw [dashed,double,thick] (L) -- (2L+2) -- (2L);
\draw [dashed,double,thick] (1) -- (2L+1) -- (L+1);
\end{tikzpicture}
\]
\end{definition}

\begin{remark} The $\mathcal{L}$-graphs $T^4_{2k}$, ${T_{2k}^4}'$ given in Figs. \ref{T2k4figure}, \ref{T2k4'figure} are of the form $\mathcal{T}_{2k}$. \end{remark}

\begin{proposition}
For $d=-7$ and for each $k$, the $\mathcal{L}$-graph ${T^4_{2k}}'$ given in Fig. \ref{T2k4'figure} is inequivalent to the $\mathcal{L}$-graph $T^4_{2k}$ given in Fig. \ref{T2k4figure}.
\end{proposition}

\begin{proof}
Let $M$,$M'$ be the matrix representatives of $T^4_{2k},{T^4_{2k}}'$ respectively. If $M$ is strongly equivalent to $M'$ then there exists a permutation matrix $P$ and a switching matrix $S$ such that \[M=PSM'S^{-1}P^{-1}\] where $S=S^{-1}=\mbox{diag}(s_1,\ldots,s_{2k})$ for $s_i\in\mathcal{L}_1=\{\pm1\}$; and there exists $\sigma\in\mathbb{S}_{2k}$ such that for matrices $X,Y$, if $X=PYP^{-1}$ then $X_{i,j}=Y_{\sigma(i),\sigma(j)}$.

Thus in general $M_{i,j}=s_{\sigma(i)}s_{\sigma(j)}M'_{\sigma(i),\sigma(j)}=\pm M'_{\sigma(i),\sigma(j)}$. Since $\omega=M_{1,2L+1}=M_{L+1,2L+1}=M_{L,2L+2}=-M_{2L,2L+2}$, considering the entries $\pm\omega$ in $M'$ we therefore require that \[\{M'_{\sigma(1),\sigma(2L+1)},M'_{\sigma(L+1),\sigma(2L+1)},M'_{\sigma(L),\sigma(2L+2)},M'_{\sigma(2L),\sigma(2L+2)}\} = \{M'_{1,2L+1},M'_{L+1,2L+1},M'_{2L+2,L},M'_{2L+2,2L}\},\] which is impossible since it implies \[\{\sigma(2L+1),\sigma(2L+2)\}=\{2L+1,L,2L\}\]

For $-M$ strongly equivalent to $M'$ we obtain the same condition, whilst for $\pm\overline{M}$ strongly equivalent to $M'$ we would require \[\{M'_{\sigma(1),\sigma(2L+1)},M'_{\sigma(L+1),\sigma(2L+1)},M'_{\sigma(L),\sigma(2L+2)},M'_{\sigma(2L),\sigma(2L+2)}\} = \{M'_{L,2L+2},M'_{2L,2L+2},M'_{2L+1,1},M'_{2L+1,L+1}\}\] which is also impossible. So $M,M'$ are necessarily inequivalent.
\end{proof}

We note the following useful computational results:

\begin{lemma}\label{T2k4left}
If $G$ is cyclotomic and induces a subgraph of the form \[
\begin{tikzpicture}[scale=1.5]
\node (1) at (0,0) [fill=black,draw,circle,label=above:{\scriptsize$1$}] {};
\node (2) at (1,0) [fill=black,draw,circle,label=above:{\scriptsize$2$}] {};
\draw (1) -- (2);
\node (L+1) at (0,-1) [fill=black,draw,circle,label=below:{\scriptsize$L+1$}] {};
\node (L+2) at (1,-1) [fill=black,draw,circle,label=below:{\scriptsize$L+2$}] {};
\draw [dotted,thick] (L+1) -- (L+2);
\draw (1) -- (L+2);
\draw [dotted,thick] (L+1) -- (2);
\node (2L+1) at (-0.5,-0.5) [fill=black,draw,circle,label=left:{\scriptsize$2L+1$}] {};
\draw [dashed,double,thick] (1) -- node[black,above left] {\scriptsize $\alpha$} (2L+1) -- node[black,below left] {\scriptsize $\beta$} (L+1);
\end{tikzpicture}
\] then $\alpha=\beta\in\mathcal{L}_2$.
\end{lemma}

\begin{lemma}\label{T2k4right}
If $G$ is cyclotomic and induces a subgraph of the form \[
\begin{tikzpicture}[scale=1.5]
\node (1) at (0,0) [fill=black,draw,circle,label=above:{\scriptsize$L-1$}] {};
\node (2) at (1,0) [fill=black,draw,circle,label=above:{\scriptsize$L$}] {};
\draw (1) -- (2);
\node (L+1) at (0,-1) [fill=black,draw,circle,label=below:{\scriptsize$2L-1$}] {};
\node (L+2) at (1,-1) [fill=black,draw,circle,label=below:{\scriptsize$2L$}] {};
\draw [dotted,thick] (L+1) -- (L+2);
\draw (1) -- (L+2);
\draw [dotted,thick] (L+1) -- (2);
\node (2L+1) at (1.5,-0.5) [fill=black,draw,circle,label=right:{\scriptsize$2L+2$}] {};
\draw [dashed,double,thick] (2) -- node[black,above right] {\scriptsize $\gamma$} (2L+1) -- node[black,below right] {\scriptsize $\delta$} (L+2);
\end{tikzpicture}
\] then $\gamma=-\delta\in\mathcal{L}_2$.
\end{lemma}

\begin{proposition}\label{T2k4twoclass}
If $G$ is a cyclotomic $\mathcal{L}$-graph of form $\mathcal{T}^4_{2k}$ then it is equivalent to the $\mathcal{L}$-graph $T^4_{2k}$ given in Fig. \ref{T2k4figure} or ($d=-7$ only) ${T^4_{2k}}'$ given in Fig. \ref{T2k4'figure}.
\end{proposition}

\begin{proof}
For $k\ge5$ the result is immediate: for the vertex numbering given in Definition \ref{T2k4def}, vertices $1,\ldots,2L$ are a cylinder of length at least 4, so by Corollary \ref{cylinderfixall} $G$ is equivalent to an $\mathcal{L}$-graph of form 
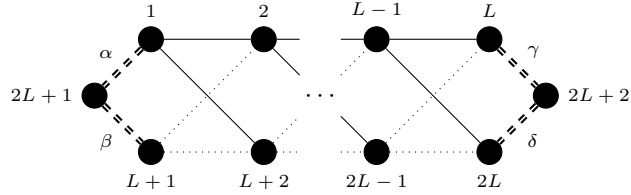
\[
\begin{tikzpicture}[scale=1.5]
\node (1) at (1,0) [fill=black,draw,circle,label=above:{\scriptsize$1$}] {};
\node (3) at (2,0) [fill=black,draw,circle,label=above:{\scriptsize$2$}] {};
\node (3a) at (2.4,0) {};
\node (3b) at (2.4,-0.4) {};
\draw (1) -- (3) -- (3a); 
\draw (3) -- (3b);
\node (L+1) at (1,-1) [fill=black,draw,circle,label=below:{\scriptsize$L+1$}] {};
\node (L+3) at (2,-1) [fill=black,draw,circle,label=below:{\scriptsize$L+2$}] {};
\node (k+3a) at (2.4,-1) {};
\node (k+3b) at (2.4,-0.6) {};
\draw [dotted] (L+1) -- (L+3) -- (k+3a);
\draw [dotted] (L+3) -- (k+3b);
\draw (1) -- (L+3);
\draw [dotted] (L+1) -- (3);
\node (spacer) at (2.5,-0.5) {$\cdots$};
\node (k-1b) at (2.6,-0.4) {};
\node (2k-1b) at (2.6,-0.6) {};
\node (k-1a) at (2.6,0) {};
\node (L-1) at (3,0) [fill=black,draw,circle,label=above:{\scriptsize$L-1$}] {};
\node (L) at (4,0) [fill=black,draw,circle,label=above:{\scriptsize$L$}] {};
\draw (k-1a) -- (L-1) -- (L);
\node (2k-1a) at (2.6,-1) {};
\node (2L-1) at (3,-1) [fill=black,draw,circle,label=below:{\scriptsize$2L-1$}] {};
\node (2L) at (4,-1) [fill=black,draw,circle,label=below:{\scriptsize$2L$}] {};
\draw [dotted] (2k-1a) -- (2L-1) -- (2L);
\draw (2k-1b) -- (2L-1) (L-1) -- (2L);
\draw [dotted] (k-1b) -- (L-1) (2L-1) -- (L); 
\node (2L+1) at (0.5,-0.5) [fill=black,draw,circle,label=left:{\scriptsize$2L+1$}] {};
\node (2L+2) at (4.5,-0.5) [fill=black,draw,circle,label=right:{\scriptsize$2L+2$}] {};
\draw [dashed,double,thick] (L) -- node[black,above right] {\scriptsize $\gamma$} (2L+2) -- node[black,below right] {\scriptsize $\delta$} (2L);
\draw [dashed,double,thick] (1) -- node[black,above left] {\scriptsize $\alpha$} (2L+1) -- node[black,below left] {\scriptsize $\beta$} (L+1);
\end{tikzpicture}
\] for some $\alpha,\beta,\gamma,\delta\in\mathcal{L}_2$ . Then by Lemmata \ref{T2k4left}, \ref{T2k4right} we have that $\alpha=\beta$ and $\gamma=-\delta$. For $d=-2$, $\mathcal{L}_2=\{\pm\sqrt{-2}\}$ so by switching at $2L+1,2L+2$ we can ensure that $\alpha=\gamma=\sqrt{-2}$, giving the $\mathcal{L}$-graph $T^4_{2k}$. 
For $d=-7$, by negation and/or conjugation $G$ is equivalent to an $\mathcal{L}$-graph with $\alpha=\omega=\frac{1}{2}+\frac{\sqrt{-7}}{2}$, and by switching at vertex $2L+2$ we can ensure $\gamma=\omega$ - giving $T^4_{2k}$ - or that $\gamma=\overline{\omega}$, giving ${T^4_{2k}}'$.

For $k=2,3$, or $4$ we can verify the result directly, after first fixing a subset of the edge labels by the equivalence operations.
\end{proof}

\subsection{$\mathcal{L}$-graphs of the form $\mathcal{C}_{2k}^{2\pm}$}
\begin{definition}
For $k\ge1$ define the $2k+1$-vertex form $\mathcal{C}_{2k}^{2\pm}$ by

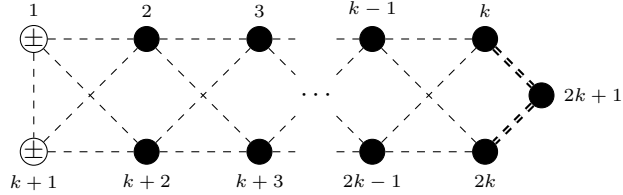
\[
\begin{tikzpicture}[scale=1.5]
\node (1) at (0,0) [fill=white,inner sep=0pt,draw,circle,label=above:{\scriptsize$1$}] {$\pm$};
\node (2) at (1,0) [fill=black,draw,circle,label=above:{\scriptsize$2$}] {};
\node (3) at (2,0) [fill=black,draw,circle,label=above:{\scriptsize$3$}] {};
\node (3a) at (2.4,0) {};
\node (3b) at (2.4,-0.4) {};
\draw [dashed] (1) -- (2) -- (3) -- (3a); 
\draw [dashed] (3) -- (3b);
\node (L+1) at (0,-1) [fill=white,inner sep=0pt,draw,circle,label=below:{\scriptsize$k+1$}] {$\pm$};
\node (L+2) at (1,-1) [fill=black,draw,circle,label=below:{\scriptsize$k+2$}] {};
\node (L+3) at (2,-1) [fill=black,draw,circle,label=below:{\scriptsize$k+3$}] {};
\node (k+3a) at (2.4,-1) {};
\node (k+3b) at (2.4,-0.6) {};
\draw [dashed] (1) -- (L+1) -- (L+2) -- (L+3) -- (k+3a);
\draw [dashed] (L+3) -- (k+3b);
\draw [dashed] (1) -- (L+2) -- (3);
\draw [dashed] (L+1) -- (2) -- (L+3);
\node (spacer) at (2.5,-0.5) {$\cdots$};
\node (k-1b) at (2.6,-0.4) {};
\node (2k-1b) at (2.6,-0.6) {};
\node (k-1a) at (2.6,0) {};
\node (L-1) at (3,0) [fill=black,draw,circle,label=above:{\scriptsize$k-1$}] {};
\node (L) at (4,0) [fill=black,draw,circle,label=above:{\scriptsize$k$}] {};
\draw [dashed] (k-1a) -- (L-1) -- (L);
\node (2k-1a) at (2.6,-1) {};
\node (2L-1) at (3,-1) [fill=black,draw,circle,label=below:{\scriptsize$2k-1$}] {};
\node (2L) at (4,-1) [fill=black,draw,circle,label=below:{\scriptsize$2k$}] {};
\draw [dashed] (2k-1a) -- (2L-1) -- (2L);
\draw [dashed] (2k-1b) -- (2L-1) -- (L);
\draw [dashed] (k-1b) -- (L-1) -- (2L);
\node (2L+1) at (4.5,-0.5) [fill=black,draw,circle,label=right:{\scriptsize$2k+1$}] {};
\draw [dashed,double,thick] (L) -- (2L+1) -- (2L);
\end{tikzpicture}
\]

\end{definition}

\begin{proposition}\label{C2k2+allequiv}
If $G$ is a cyclotomic charged $\mathcal{L}$-graph of form $\mathcal{C}^{2\pm}_{2k}$ then it is equivalent to the charged $\mathcal{L}$-graph $C^{2+}_{2k}$ given in Fig \ref{C2k+figure}.
\end{proposition}

\begin{proof}
For $k\ge5$, the result is immediate. By Lemma \ref{samesign} we have that the charges on vertices $1,k+1$ are equal; negating if necessary $G$ is equivalent to an $\mathcal{L}$-graph with both charges $+1$. Then vertices $2,\ldots,k,k+2,\ldots 2k$ are a cylinder of length at least 4, so by Corollary \ref{cylinderfixall} and switching at $1,k+1$ $G$ is equivalent to an $\mathcal{L}$-graph with edges specified as follows: 
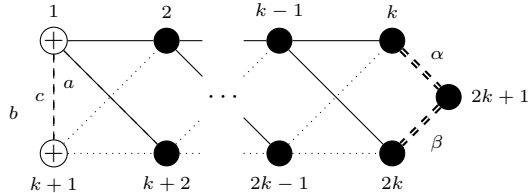
\[
\begin{tikzpicture}[scale=1.5]
\node (1) at (1,0) [fill=white,inner sep=0pt,draw,circle,label=above:{\scriptsize$1$}] {$+$};
\node (3) at (2,0) [fill=black,draw,circle,label=above:{\scriptsize$2$}] {};
\node (3a) at (2.4,0) {};
\node (3b) at (2.4,-0.4) {};
\draw (1) -- (3) -- (3a); 
\draw (3) -- (3b);
\node (L+1) at (1,-1) [fill=white,inner sep=0pt,draw,circle,label=below:{\scriptsize$k+1$}] {$+$};
\node (L+3) at (2,-1) [fill=black,draw,circle,label=below:{\scriptsize$k+2$}] {};
\node (k+3a) at (2.4,-1) {};
\node (k+3b) at (2.4,-0.6) {};
\draw [dotted] (L+1) -- (L+3) -- (k+3a);
\draw [dotted] (L+3) -- (k+3b);
\draw (1) -- (L+3);
\draw [dotted] (L+1) -- (3);
\node (spacer) at (2.5,-0.5) {$\cdots$};
\node (k-1b) at (2.6,-0.4) {};
\node (2k-1b) at (2.6,-0.6) {};
\node (k-1a) at (2.6,0) {};
\node (L-1) at (3,0) [fill=black,draw,circle,label=above:{\scriptsize$k-1$}] {};
\node (L) at (4,0) [fill=black,draw,circle,label=above:{\scriptsize$k$}] {};
\draw (k-1a) -- (L-1) -- (L);
\node (2k-1a) at (2.6,-1) {};
\node (2L-1) at (3,-1) [fill=black,draw,circle,label=below:{\scriptsize$2k-1$}] {};
\node (2L) at (4,-1) [fill=black,draw,circle,label=below:{\scriptsize$2k$}] {};
\draw [dotted] (2k-1a) -- (2L-1) -- (2L);
\draw (2k-1b) -- (2L-1) (L-1) -- (2L);
\draw [dotted] (k-1b) -- (L-1) (2L-1) -- (L); 
\node (2L+2) at (4.5,-0.5) [fill=black,draw,circle,label=right:{\scriptsize$2k+1$}] {};
\draw [dashed,double,thick] (L) -- node[black,above right] {\scriptsize $\alpha$} (2L+2) -- node[black,below right] {\scriptsize $\beta$} (2L);
\draw [dashed] (1) -- node[left,black] {\scriptsize$c$} (L+1) (1) -- node[above left,black,xshift=-1em] {\scriptsize$a$} (L+3) (L+1) -- node[below left,black,xshift=-1em] {\scriptsize$b$} (2);
\end{tikzpicture}
\]

for some $a,b,c\in\mathcal{L}_1$, $\alpha,\beta\in\mathcal{L}_2$. But, by Corollary \ref{C2k++cor}, the subgraph induced on vertices $1,2,3,4,k+1,k+2,k+3,k+4$ is cyclotomic if and only if $a=c=-b=1$. By complex conjugation and/or switching at $2k+1$, we can ensure $\alpha=\sqrt{-2}$ or $\frac{1}{2}+\frac{\sqrt{-7}}{2}$ for $d=-2,-7$ respectively; by Lemma \ref{T2k4right}, $\beta=-\alpha$. Thus we recover the charged $\mathcal{L}$-graph $C^{2+}_{2k}$ as claimed.

If $k=1,2,3$ or $4$ then, by Lemma \ref{samesign} and fixing a subset of the edge labels under equivalence, the result can be verified directly.

\end{proof}

\section{Classification of $4$-cyclotomic $\mathcal{L}$-graphs up to form}

We complete the proof of Theorems \ref{cyc_class_neg2}, \ref{cyc_class_neg7} by demonstrating the following:

\begin{proposition}\label{4cyc_class}
If $G$ is a $4$-cyclotomic $\mathcal{L}$-graph with at least one edge label from $\mathcal{L}_2$ and all such edges occuring in isolated pairs, then $G$ is of the form $\mathcal{T}_{2k}^4$ or $\mathcal{C}_{2k}^{2\pm}$.
\end{proposition}

\begin{definition}
If we have a cyclotomic $\mathcal{L}$-graph G on vertices $v_1\ldots v_k$, we will describe the extension of G by vertices $x_1\ldots x_n$ and corresponding edges as a \emph{saturating extension} if all vertices $v_1\ldots v_k$ then have weighted degree four; the $x_i$ needn't also be saturated.
\end{definition}
Trivially, any subgraph $G'$ of a $4$-cyclotomic $\mathcal{L}$-graph $G$ can be grown to $G$ by a saturating extension- simply reintroduce all missing vertices and edges. We thus describe a saturating extension by $x_1\ldots x_n$ as \emph{minimal} if omitting any one of the $x_i$ and its corresponding edges gives a non-saturating extension (that is, each $x_i$ is necessary to saturate some $v_j$). Note that a minimal saturating extension corresponds to some sequence of saturating additions. 

\begin{proposition}\label{cansat}
 Any 4-cyclotomic $\mathcal{L}$-graph $G$ can be grown from any of its induced subgraphs by a sequence of minimal saturating extensions.\end{proposition}

\begin{proposition}[Base Step]\label{ucthm2a}
Given the graph $\begin{tikzpicture}
\node (a) at (0,0) [fill=black,draw,circle] {};
\node (b) at (1,0) [fill=black,draw,circle] {};
\node (c) at (2,0) [fill=black,draw,circle] {};
\draw [double,dashed,thick] (a) -- (b) -- (c);
\end{tikzpicture}$, the only possible minimal saturating extensions are maximal graphs of the form $\mathcal{T}_{4}^4$:
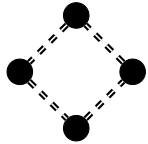
\[
\begin{tikzpicture}[scale=1.5]
\node (a) at (0,0) [fill=black,draw,circle] {};
\node (b) at (-0.5,-0.5) [fill=black,draw,circle] {};
\node (c) at (0.5,-0.5) [fill=black,draw,circle] {};
\node (d) at (0,-1) [fill=black,draw,circle] {};
\draw [double,dashed,thick] (a) -- (b) -- (d) -- (c) -- (a);
\end{tikzpicture}
\]
or maximal graphs of the form $\mathcal{C}_4^{2\pm}$:
\[
\begin{tikzpicture}[scale=1.5]
\node (L-1) at (3,0) [fill=white,inner sep=0pt,draw,circle] {$\pm$};
\node (L) at (4,0) [fill=black,draw,circle] {};
\draw [dashed] (L-1) -- (L);
\node (2L-1) at (3,-1) [fill=white,inner sep=0pt,draw,circle] {$\pm$};
\node (2L) at (4,-1) [fill=black,draw,circle] {};
\draw [dashed] (L) -- (2L-1) -- (2L) -- (L-1) -- (2L-1);
\node (2L+2) at (4.5,-0.5) [fill=black,draw,circle] {};
\draw [dashed,double,thick] (L) -- (2L+2) -- (2L);
\end{tikzpicture}
\]
or nonmaximal \emph{chains of length one}: 
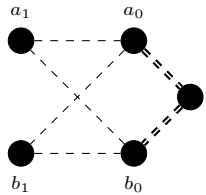
\[
\begin{tikzpicture}[scale=1.5]
\node (L-1) at (3,0) [fill=black,draw,circle,label=above:{\scriptsize$a_1$}] {};
\node (L) at (4,0) [fill=black,draw,circle,label=above:{\scriptsize$a_0$}] {};
\draw [dashed] (L-1) -- (L);
\node (2L-1) at (3,-1) [fill=black,draw,circle,label=below:{\scriptsize$b_1$}] {};
\node (2L) at (4,-1) [fill=black,draw,circle,label=below:{\scriptsize$b_0$}] {};
\draw [dashed] (L) -- (2L-1) -- (2L) -- (L-1);
\node (2L+2) at (4.5,-0.5) [fill=black,draw,circle] {};
\draw [dashed,double,thick] (L) -- (2L+2) -- (2L);
\end{tikzpicture}
\]
\end{proposition}

\begin{proof} Let the extension set be $x_1,\ldots,x_n$. If some $x_i$ is joined to $a_0$ or $b_0$ by an edge of weight 2 then $x_i$ is uncharged by Lemma \ref{nonMaxSubGraphsw2} (A) and so we have a path of three consecutive weight 2 edges, forcing (by part (D) of the same Lemma) the graph to be of form $\mathcal{T}_4^4$ as required.

Thus we may assume each edge from an $x_j$ to $a_0,b_0$ is of weight 1; to satisfy both minimality and saturation this forces $n=2,3$ or $4$. However, if $n\neq2$ then there exists a neighbour of $a_0$ which is not a neighbour of $b_0$, which induces a subgraph of the form 
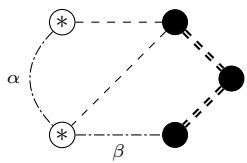
\[
\begin{tikzpicture}[scale=1.5]
\node (a0) at (0,0) [fill=black,draw,circle] {};
\node (z) at (0.5,-0.5) [fill=black,draw,circle] {};
\node (b0) at (0,-1) [fill=black,draw,circle] {};
\node (x1) at (-1,0) [fill=white,inner sep=1pt,draw,circle] {$*$};
\node (x2) at (-1,-1) [fill=white,inner sep=1pt,draw,circle] {$*$};
\draw [double,dashed,thick] (a0) -- (z) -- (b0);
\draw [dashed] (x1) -- (a0) -- (x2);
\draw [dashdot] (x2) -- node [below,black] {\scriptsize$\beta$} (b0); 
\draw [dashdot] (x1) .. controls +(225:0.5) and +(135:0.5) .. node [left,black] {\scriptsize$\alpha$} (x2); 
\end{tikzpicture}
\] However, no such $\mathcal{L}$-graph is cyclotomic for $\alpha\in\mathcal{L}$, $\beta\in\mathcal{L}_1\cup\{0\}$.
\end{proof}

\begin{proposition}[Inductive Step]\label{ucthm2b}
Given a \emph{chain of length k}:
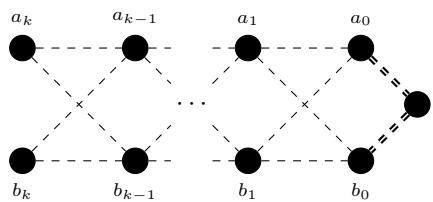
\[
\begin{tikzpicture}[scale=1.5]
\node (2) at (1,0) [fill=black,draw,circle,label=above:{\scriptsize$a_k$}] {};
\node (3) at (2,0) [fill=black,draw,circle,label=above:{\scriptsize$a_{k-1}$}] {};
\node (3a) at (2.4,0) {};
\node (3b) at (2.4,-0.4) {};
\draw [dashed] (2) -- (3) -- (3a); 
\draw [dashed] (3) -- (3b);
\node (L+2) at (1,-1) [fill=black,draw,circle,label=below:{\scriptsize$b_k$}] {};
\node (L+3) at (2,-1) [fill=black,draw,circle,label=below:{\scriptsize$b_{k-1}$}] {};
\node (k+3a) at (2.4,-1) {};
\node (k+3b) at (2.4,-0.6) {};
\draw [dashed] (L+2) -- (L+3) -- (k+3a);
\draw [dashed] (L+3) -- (k+3b);
\draw [dashed] (L+2) -- (3);
\draw [dashed] (2) -- (L+3);
\node (spacer) at (2.5,-0.5) {$\cdots$};
\node (k-1b) at (2.6,-0.4) {};
\node (2k-1b) at (2.6,-0.6) {};
\node (k-1a) at (2.6,0) {};
\node (L-1) at (3,0) [fill=black,draw,circle,label=above:{\scriptsize$a_1$}] {};
\node (L) at (4,0) [fill=black,draw,circle,label=above:{\scriptsize$a_0$}] {};
\draw [dashed] (k-1a) -- (L-1) -- (L);
\node (2k-1a) at (2.6,-1) {};
\node (2L-1) at (3,-1) [fill=black,draw,circle,label=below:{\scriptsize$b_1$}] {};
\node (2L) at (4,-1) [fill=black,draw,circle,label=below:{\scriptsize$b_0$}] {};
\draw [dashed] (2k-1a) -- (2L-1) -- (2L);
\draw [dashed] (2k-1b) -- (2L-1) -- (L);
\draw [dashed] (k-1b) -- (L-1) -- (2L);
\node (2L+2) at (4.5,-0.5) [fill=black,draw,circle] {};
\draw [dashed,double,thick] (L) -- (2L+2) -- (2L);
\end{tikzpicture}
\]
the only possible minimal saturating extensions are maximal graphs of the form $\mathcal{T}_{2(k+2)}^4$:
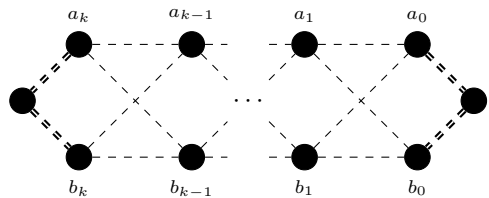
\[
\begin{tikzpicture}[scale=1.5]
\node (2) at (1,0) [fill=black,draw,circle,label=above:{\scriptsize$a_k$}] {};
\node (3) at (2,0) [fill=black,draw,circle,label=above:{\scriptsize$a_{k-1}$}] {};
\node (3a) at (2.4,0) {};
\node (3b) at (2.4,-0.4) {};
\draw [dashed] (2) -- (3) -- (3a); 
\draw [dashed] (3) -- (3b);
\node (L+2) at (1,-1) [fill=black,draw,circle,label=below:{\scriptsize$b_k$}] {};
\node (L+3) at (2,-1) [fill=black,draw,circle,label=below:{\scriptsize$b_{k-1}$}] {};
\node (k+3a) at (2.4,-1) {};
\node (k+3b) at (2.4,-0.6) {};
\draw [dashed] (L+2) -- (L+3) -- (k+3a);
\draw [dashed] (L+3) -- (k+3b);
\draw [dashed] (L+2) -- (3);
\draw [dashed] (2) -- (L+3);
\node (spacer) at (2.5,-0.5) {$\cdots$};
\node (k-1b) at (2.6,-0.4) {};
\node (2k-1b) at (2.6,-0.6) {};
\node (k-1a) at (2.6,0) {};
\node (L-1) at (3,0) [fill=black,draw,circle,label=above:{\scriptsize$a_1$}] {};
\node (L) at (4,0) [fill=black,draw,circle,label=above:{\scriptsize$a_0$}] {};
\draw [dashed] (k-1a) -- (L-1) -- (L);
\node (2k-1a) at (2.6,-1) {};
\node (2L-1) at (3,-1) [fill=black,draw,circle,label=below:{\scriptsize$b_1$}] {};
\node (2L) at (4,-1) [fill=black,draw,circle,label=below:{\scriptsize$b_0$}] {};
\draw [dashed] (2k-1a) -- (2L-1) -- (2L);
\draw [dashed] (2k-1b) -- (2L-1) -- (L);
\draw [dashed] (k-1b) -- (L-1) -- (2L);
\node (2L+2) at (4.5,-0.5) [fill=black,draw,circle] {};
\draw [dashed,double,thick] (L) -- (2L+2) -- (2L);
\node (2L+1) at (0.5,-0.5) [fill=black,draw,circle] {};
\draw [dashed,double,thick] (2) -- (2L+1) -- (L+2);
\end{tikzpicture}
\]
or maximal graphs of the form $\mathcal{C}_{2(k+2)}^{2\pm}$:
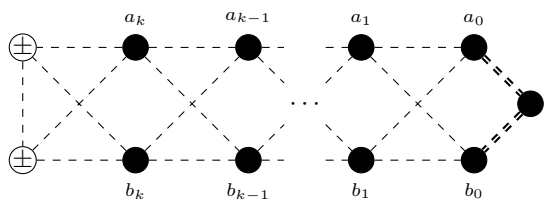
\[
\begin{tikzpicture}[scale=1.5]
\node (2) at (1,0) [fill=black,draw,circle,label=above:{\scriptsize$a_k$}] {};
\node (3) at (2,0) [fill=black,draw,circle,label=above:{\scriptsize$a_{k-1}$}] {};
\node (3a) at (2.4,0) {};
\node (3b) at (2.4,-0.4) {};
\draw [dashed] (2) -- (3) -- (3a); 
\draw [dashed] (3) -- (3b);
\node (L+2) at (1,-1) [fill=black,draw,circle,label=below:{\scriptsize$b_k$}] {};
\node (L+3) at (2,-1) [fill=black,draw,circle,label=below:{\scriptsize$b_{k-1}$}] {};
\node (k+3a) at (2.4,-1) {};
\node (k+3b) at (2.4,-0.6) {};
\draw [dashed] (L+2) -- (L+3) -- (k+3a);
\draw [dashed] (L+3) -- (k+3b);
\draw [dashed] (L+2) -- (3);
\draw [dashed] (2) -- (L+3);
\node (spacer) at (2.5,-0.5) {$\cdots$};
\node (k-1b) at (2.6,-0.4) {};
\node (2k-1b) at (2.6,-0.6) {};
\node (k-1a) at (2.6,0) {};
\node (L-1) at (3,0) [fill=black,draw,circle,label=above:{\scriptsize$a_1$}] {};
\node (L) at (4,0) [fill=black,draw,circle,label=above:{\scriptsize$a_0$}] {};
\draw [dashed] (k-1a) -- (L-1) -- (L);
\node (2k-1a) at (2.6,-1) {};
\node (2L-1) at (3,-1) [fill=black,draw,circle,label=below:{\scriptsize$b_1$}] {};
\node (2L) at (4,-1) [fill=black,draw,circle,label=below:{\scriptsize$b_0$}] {};
\draw [dashed] (2k-1a) -- (2L-1) -- (2L);
\draw [dashed] (2k-1b) -- (2L-1) -- (L);
\draw [dashed] (k-1b) -- (L-1) -- (2L);
\node (2L+2) at (4.5,-0.5) [fill=black,draw,circle] {};
\draw [dashed,double,thick] (L) -- (2L+2) -- (2L);
\node (x1) at (0,0) [fill=white,draw,circle,inner sep=0pt] {$\pm$};
\node (x2) at (0,-1) [fill=white,draw,circle,inner sep=0pt] {$\pm$};
\draw [dashed] (x2) -- (2) -- (x1) -- (x2) -- (L+2) -- (x1);
\end{tikzpicture}
\]
or a nonmaximal \emph{chain of length $k+1$}.  
\end{proposition}

\begin{proof}
Let $X=\{x_1\,\ldots x_n\}$ be the saturating set.

For $k=1$, we note that there are no cyclotomic graphs of the form 
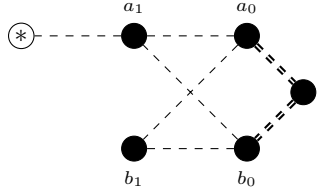
\[
\begin{tikzpicture}[scale=1.5]
\node (a0) at (0,0) [fill=black,draw,circle,label=above:{\scriptsize$a_0$}] {};
\node (z) at (0.5,-0.5) [fill=black,draw,circle] {};
\node (b0) at (0,-1) [fill=black,draw,circle,label=below:{\scriptsize$b_0$}] {};
\node (x1) at (-1,0) [fill=black,draw,circle,label=above:{\scriptsize$a_1$}] {};
\node (x2) at (-1,-1) [fill=black,draw,circle,label=below:{\scriptsize$b_1$}] {};
\node (x) at (-2,0) [fill=white,inner sep=1pt,draw,circle] {$*$};
\draw [double,dashed,thick] (a0) -- (z) -- (b0);
\draw [dashed] (x1) -- (a0) -- (x2) -- (b0) -- (x1) -- (x);
\end{tikzpicture}
\]  
Such a graph is necessarily induced if $n=3,4$ or $n=2$ with an edge of weight 2 between one of the $x_i$ and either $a_1$ or $b_1$. Thus we either have $n=1$ which forces a graph of form $\mathcal{T}_6^4$, or $n=2$ with all new edges of weight 1. Such a graph is cyclotomic only if it's a maximal graph of form $\mathcal{C}_6^{2+}$ or a chain of length 2, as required.

Otherwise $k\ge2$ and we note the following result:

\begin{lemma}\label{weight2G} The $\mathcal{L}$-graph $H$
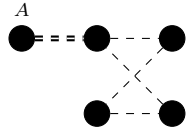
\[
\begin{tikzpicture}
\node (2) at (1,0) [fill=black,draw,circle,label=above:{\scriptsize$A$}] {};
\node (3) at (2,0) [fill=black,draw,circle] {};
\node (4) at (3,0) [fill=black,draw,circle] {};
\node (5) at (2,-1) [fill=black,draw,circle] {};
\node (6) at (3,-1) [fill=black,draw,circle] {};
\draw [dashed] (3) -- (4) -- (5) -- (6) -- (3);
\draw [dashed,double,thick] (2) -- (3);
\end{tikzpicture}
\]
cannot be an induced subgraph of a 4-cyclotomic $\mathcal{L}$-graph $G$ whose weight 2 edges arise as isolated pairs.
\end{lemma}
\begin{proof} Vertex $A$ would necessarily have additional neighbours in $G$, but the only cyclotomic possibilities induce an isolated weight 2 edge.\end{proof}

By Lemma \ref{weight2G} if $n\ge1$ we may assume no edges of weight 2 join $a_k,b_k$ to any of the $x_i$. But then Corollary \ref{kgt1nocyc} allows us to exclude $n=3$ or $4$ since there would be a neighbour of $a_k$ not neighbouring $b_k$. So we either have $n=1$, which to ensure saturation forces a graph of form $\mathcal{T}_{2(k+2)}^4$, or $n=2$ with a graph of form 
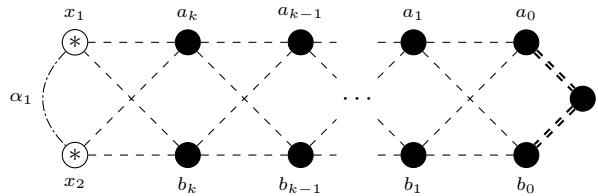
\[
\begin{tikzpicture}[scale=1.5]
\node (2) at (1,0) [fill=black,draw,circle,label=above:{\scriptsize$a_k$}] {};
\node (3) at (2,0) [fill=black,draw,circle,label=above:{\scriptsize$a_{k-1}$}] {};
\node (3a) at (2.4,0) {};
\node (3b) at (2.4,-0.4) {};
\draw [dashed] (2) -- (3) -- (3a); 
\draw [dashed] (3) -- (3b);
\node (L+2) at (1,-1) [fill=black,draw,circle,label=below:{\scriptsize$b_k$}] {};
\node (L+3) at (2,-1) [fill=black,draw,circle,label=below:{\scriptsize$b_{k-1}$}] {};
\node (k+3a) at (2.4,-1) {};
\node (k+3b) at (2.4,-0.6) {};
\draw [dashed] (L+2) -- (L+3) -- (k+3a);
\draw [dashed] (L+3) -- (k+3b);
\draw [dashed] (L+2) -- (3);
\draw [dashed] (2) -- (L+3);
\node (spacer) at (2.5,-0.5) {$\cdots$};
\node (k-1b) at (2.6,-0.4) {};
\node (2k-1b) at (2.6,-0.6) {};
\node (k-1a) at (2.6,0) {};
\node (L-1) at (3,0) [fill=black,draw,circle,label=above:{\scriptsize$a_1$}] {};
\node (L) at (4,0) [fill=black,draw,circle,label=above:{\scriptsize$a_0$}] {};
\draw [dashed] (k-1a) -- (L-1) -- (L);
\node (2k-1a) at (2.6,-1) {};
\node (2L-1) at (3,-1) [fill=black,draw,circle,label=below:{\scriptsize$b_1$}] {};
\node (2L) at (4,-1) [fill=black,draw,circle,label=below:{\scriptsize$b_0$}] {};
\draw [dashed] (2k-1a) -- (2L-1) -- (2L);
\draw [dashed] (2k-1b) -- (2L-1) -- (L);
\draw [dashed] (k-1b) -- (L-1) -- (2L);
\node (2L+2) at (4.5,-0.5) [fill=black,draw,circle] {};
\draw [dashed,double,thick] (L) -- (2L+2) -- (2L);
\node (x1) at (0,0) [fill=white,draw,circle,inner sep=1pt,label=above:{\scriptsize$x_1$}] {$*$};
\node (x2) at (0,-1) [fill=white,draw,circle,inner sep=1pt,label=below:{\scriptsize$x_2$}] {$*$};
\draw [dashed] (x2) -- (2) -- (x1) (x2) -- (L+2) -- (x1);
\draw [dashdot] (x1) .. controls +(225:0.5) and +(135:0.5) .. node [left,black] {\scriptsize$\alpha_1$} (x2); 
\end{tikzpicture}
\]
By interlacing, it suffices to check the possible subgraphs on vertices $x_1, a_k, a_{k-1},x_2, b_k, b_{k-1}$ for cyclotomicity, which confirms that the only possibilities are a graph of form $\mathcal{C}_{2(k+2)}^{2\pm}$ or a chain of length $k+1$, as required.
\end{proof}

Thus Proposition \ref{4cyc_class} holds: $G$ is either of form $C_2^{2+}$, or it induces a subgraph of form $\begin{tikzpicture}
\node (a) at (0,0) [fill=black,draw,circle] {};
\node (b) at (1,0) [fill=black,draw,circle] {};
\node (c) at (2,0) [fill=black,draw,circle] {};
\draw [double,dashed,thick] (a) -- (b) -- (c);
\end{tikzpicture}$. By Proposition \ref{cansat} it can therefore be grown by a sequence of minimal saturating extensions, terminating with $G$, which is maximal. Since a chain is not maximal, by Proposition \ref{ucthm2a} and Theorem \ref{ucthm2b} $G$ must be of form $\mathcal{T}_{2k}$ or $\mathcal{C}_{2k}^{2\pm}$ for some $k$.

So we have completed the proof of Theorems \ref{cyc_class_neg2}, \ref{cyc_class_neg7}: any maximal cyclotomic $\mathcal{L}$-graph $G$ for $d=-2,-7$ is a charged signed graph unless it has an edge label from $\mathcal{L}_2\cup\mathcal{L}_3\cup\mathcal{L}_4$; edges from $L_4$ or $L_3$ force $G$ to be equivalent to one of $S_2, S_2^*, S_2', S_4'$ by the results of Sections \ref{L4section}, \ref{L3section}; if there are any isolated weight 2 edges then by Lemmata \ref{maxis4cyc_C}, \ref{maxis4cyc_B} $G$ is equivalent to $S_4^*, S_6^\dag$ or $S_8^*$, whilst if there are three or more consecutive weight 2 edges then Lemma \ref{maxis4cyc_E} $G$ is equivalent to $T_4^4$ or ${T_4^4}'$; otherwise all weight 2 edges occur in isolated pairs, so Theorem \ref{4cyc_class} applies and by Theorems \ref{T2k4twoclass}, \ref{C2k2+allequiv} such a graph is equivalent to $T_{2k}^4$, ${T_{2k}^4}'$ or $C_{2k}^{2+}$ for an appropriate $k$.

\section{Existence of Maximal Supergraphs}
In this Section we prove Theorem \ref{containedinmax}. 

Let $G$ be a connected cyclotomic $\mathcal{L}$-graph. If $G$ is a charged signed graph, then the Theorem holds by Theorem \ref{McSmThm1} or \ref{McSmThm2}. If $G$ contains an edge of weight three or four then the results of Section \ref{w34edgeclassification} suffice; Theorem \ref{containedinmax} therefore holds for $d=-11$ or $d=-15$ and for $d=-2$ and $-7$ we may restrict our attention to $G$ a nonmaximal $\mathcal{L}$-graph containing edges of weight at most two, with at least one such edge.

If all vertices of $G$ have weight four then it is maximal and we are done; otherwise, by the results of Section \ref{section6}, $G$ admits a cyclotomic extension. Repeating this process, we either generate a $4$-cyclotomic supergraph of $G$ as desired, or a cyclotomic supergraph $G^*$ with at least eight vertices, at least one of which has weight less than four (If $G$ has eight or more vertices, take $G^*=G$). Since $G$ contains at least one edge of weight 2, so does $G*$, joining vertices $u,v$. By the Gram vector constructions we may extend $G^*$ further to ensure that $u$ and $v$ have weight four: since $G^*$ has at least eight vertices the subgraph exclusion results of Lemma \ref{nonMaxSubGraphsw2} force $u,v$ to be contained in an isolated pair of weight two edges. 

Thus $G^*$ induces a subgraph $H$ of the form $\begin{tikzpicture} \node (a) at (0,0) [fill=black,draw,shape=circle] {}; \node (b) at (0.9,0) [fill=black,draw,shape=circle] {}; \node (c) at (1.8,0) [fill=black,draw,shape=circle] {}; \draw [dashed,double,thick] (a) -- (b) -- (c); \end{tikzpicture}$. Let the vertices of $G^*$ be $x_1,\ldots,x_n$: by the results of Section \ref{section6} there exists a finite $m$ and vertices $y_1,\ldots y_m$ such that the $y_j$ saturate the $x_i$. Thus the graph $G^\dag$ on vertices $x_1,\ldots x_n,y_1,\ldots,y_m$ is a saturating extension of $G^*$. Thus it is also a saturating extension of $H$, and so $G^\dag$ can be recovered from $H$ by a sequence of minimal saturating extensions. But by Propositions \ref{ucthm2a}, \ref{ucthm2b} this forces $G^\dag$ to be of form $\mathcal{T}_{2k}^4$, $\mathcal{C}_{2k}^{2\pm}$ or a chain of length $k$ for some $k$. $G^*$ is therefore either $4$-cyclotomic or (if a chain) contained in a finite maximal 4-cyclotomic $\mathcal{L}$-graph: since $G^\dag$ is a supergraph of $G^*$ which contained $G$, we are done.

\begin{spacing}{0.4}

\end{spacing}
\end{document}